\documentstyle[12pt,amscd]{amsart}
\newtheorem{pr}{Proposition}[subsection]
\newtheorem{lm}{Lemma}[subsection]
\newtheorem{tm}{Theorem}[subsection]
\newcommand{\proj}{\bold P}
\newcommand{\grass}{\bold G}
\newcommand{\barr}{\overline}

\newcommand{\rarr}{\rightarrow}
\newcommand{\oh}{{\cal{O}}}
\newcommand{\mgb}{\barr{M_g}}
\newcommand{\Q}{Q_g(\mu,n,f)}
\newcommand{\tat}{{\theta}_{\bold P}}
\newcommand{\al}{\alpha}
\newcommand{\bay}{\beta}

\newcommand{\om}{\omega}
\newcommand{\sumo}{\bigoplus}
\newcommand{\paren}[1]{\left( #1 \right)}
\newcommand{\tp}{\bold{C}^n \otimes}
\newcommand{\E}{\cal{E}}
\newcommand{\K}{\cal{K}}
\newcommand{\F}{\cal{F}}

\begin{document}
\title{A Compactification over $\mgb$ of the Universal Moduli Space
of Slope-Semistable Vector Bundles. \\}
\author{Rahul Pandharipande$^1$}
\date{8 December 1994}
\maketitle
\pagestyle{plain}
\baselineskip=16pt
\footnotetext[1]{Supported by an NSF Graduate Fellowship. }

\baselineskip=11pt
\setcounter{section}{-1}
\baselineskip=16pt
\section{Introduction}
\subsection{Compactifications of Moduli Problems}
\label{cmp}
Initial statements of moduli problems in algebraic geometry often
do not yield compact moduli spaces. For example, the moduli space
$M_g$ of nonsingular, genus $g\geq 2$ curves is open. Compact
moduli spaces are desired for several reasons. Degeneration arguments
in moduli require compact spaces. Also, there are more techniques
available to study the
 global geometry of compact
spaces. It is therefore valuable to find
natural compactifications of open moduli problems. In the case of
$M_g$, there is a remarkable compactification due to P. Deligne and D. Mumford.
A connected, reduced, nodal curve $C$ of arithmetic genus $g\geq 2$
is {\em Deligne-Mumford stable}
if each nonsingular rational component contains at least three nodes of $C$.
$\barr{M_g}$, the moduli space of Deligne-Mumford stable genus $g$ curves, is
compact and includes $M_g$ as a dense open set.

There is a natural notion of stability for a vector bundle $E$ on a nonsingular
curve $C$.
Let the slope $\mu$ be defined as follows:
$\mu(E)=degree(E)/rank(E)$.
$E$ is {\em slope-stable (slope-semistable)} if
\begin{equation}
\label{ahab}
\mu(F)< (\leq)\  \mu(E)
\end{equation}
for every proper subbundle
$F$ of $E$.
When the degree and rank are not coprime, the moduli space of slope-stable
bundles is open. $U_{C}(e,r)$, the moduli space of slope-semistable bundles of
degree
$e$ and rank $r$, is compact. An open set of $U_{C}(e,r)$ corresponds
bijectively
to isomorphism classes of stable bundles. In general, points of $U_C(e,r)$
correspond
to equivalence classes (see section (\ref{sfor}) ) of semistable bundles.

The moduli problem of pairs $(C,E)$ where $E$ is a slope-semistable
vector bundle on a
 nonsingular curve  $C$ can not be compact. No allowance is made for curves
that degenerate to nodal curves.  A natural compactification of this
moduli problem of pairs is presented here.

\subsection{Compactification of the Moduli Problem of Pairs}
Let $\bold{C}$ be a fixed algebraically closed field. As before,
Let $M_g$ be the moduli space of nonsingular, complete, irreducible,
genus $g\geq2$ curves over the field $\bold{C}$.
For each $[C]\in M_g$, there
is a natural  projective variety, $U_C(e,r)$, parametrizing
degree $e$,
rank $r$,
slope-semistable vector bundles (up to equivalence) on $C$.
For $g\geq2$,
let $U_g(e,r)$ be the set of equivalence classes of  pairs $(C,E)$
where $[C] \in M_g$
and $E$ is a slope-semistable vector bundle on $C$
of the specified degree and rank.
A good compactification, K, of the moduli set of pairs $U_g(e,r)$
should satisfy at least the following conditions:
\begin{enumerate}
\item [(i.)] $K$ is a projective variety that functorially
parametrizes a class of geometric objects.
\item [(ii.)] $U_g(e,r)$ functorially corresponds to an open
dense subset of $K$.
\item [(iii.)] There exists a morphism $\eta: K \rarr \mgb$
such that the natural diagram commutes:
\begin{equation*}
\begin{CD}
U_g(e,r) @>>> K\\
@VVV @VV{\eta}V \\
M_g @>>> \mgb
\end{CD}
\end{equation*}
\item [(iv.)] For each $[C]\in M_g$, there exists a functorial
isomorphism $$\eta^{-1}([C]) \cong U_C(e,r)/Aut(C).$$
\end{enumerate}
The main result of this paper is the construction of a
projective variety $\barr{U_g(e,r)}$ that parametrizes
equivalence classes of slope-semistable, torsion free
sheaves on Deligne-Mumford stable, genus $g$ curves and
satisfies conditions (i-iv) above.

The definition of slope-semistability of torsion free sheaves
(due to C. Seshadri) is given in section (\ref{sfor}).

\subsection{The Method of Construction}
\label{moc}
An often successful approach to moduli constructions in algebraic
geometry involves two steps. In the first step, extra data
is added to rigidify the moduli problem. With the additional
data, the new moduli problem is solved by a  Hilbert or Quot
scheme. In the second step, the extra data is removed by a group
quotient. Geometric Invariant Theory is used to study the
quotient problem in the category of algebraic schemes. In good cases,
the final quotient is the desired moduli space.

In order to rigidify the moduli problem of pairs, the following
data is added to $(C,E)$:
\begin{enumerate}
\item[(i.)] An isomorphism
$\bold{C}^{N+1} \stackrel {\sim}{\rarr} H^0(C, \om_C^{10})$,
\item[(ii.)] An isomorphism
$\bold{C}^n \stackrel{\sim}{\rarr} H^0(C,E)$.
\end{enumerate}
Note $\om_C$ is the canonical bundle of $C$.
The numerical invariants of the moduli problem of pairs are
the genus $g$, degree $e$, and rank $r$. The rigidified
problem should have no more numerical invariants. Hence,
$N$ and $n$  must be determined by $g$, $e$, and $r$.
Certainly, $N=10(2g-2)-g$ by Riemann-Roch. It is assumed
$H^1(E)=0$ and $E$ is generated by global sections.
In the end, it is checked  these assumptions are consequences of the
stability condition for sufficiently high degree bundles.
We see $n=\chi(E)= e+r(1-g)$.

The isomorphism
of (i) canonically embeds $C$ in $\proj^{N}=\proj^{N}_{\bold C}$.
The isomorphism of (ii) exhibits $E$ as a canonical quotient
$$ \bold{C}^n \otimes \oh_C \rarr E \rarr 0.$$
The basic parameter spaces in algebraic geometry are Hilbert
and Quot schemes.
Subschemes of a fixed scheme $X$ are parametrized by Hilbert schemes
$Hilb(X)$.
Quotients of a fixed sheaf $F$ on $X$ are parametrized by Quot schemes
$Quot(X,F)$.
The rigidified curve $C$ can be parametrized by a Hilbert
scheme $H$ of curves in $\proj^N$, and the rigidified bundle $E$ can
be parametrized by a Quot scheme $Quot(C,\bold{C}^n\otimes \oh_C)$
of quotients on $C$.
In fact, Quot schemes can be defined in a relative context. Let $U_H$ be
the universal curve over the Hilbert scheme $H$.
The family of Quot schemes, $Quot(C,\bold{C}^n\otimes \oh_C)$,
defined as $C\hookrightarrow \proj^N$
varies in $H$ is simply the relative
Quot scheme of the universal curve over the Hilbert scheme:
$Quot(U_H\rarr H, \bold{C}^n\otimes \oh_{U_H})$. This relative
Quot scheme is the parameter space of the rigidified pairs (up to
scalars).
The Quot scheme set up is discussed in detail in section (\ref{rqs}).

The actions of $GL_{N+1}(\bold{C})$ on $\bold{C}^{N+1}$ and $GL_{n}(\bold{C})$
on $\bold{C}^n$ yield an action of $GL_{N+1}\times GL_{n}$ on the
rigidified data.
There is an induced product action on the relative Quot scheme.
It is easily seen the scalar elements of the groups act
trivially on the Quot scheme.  $\barr{U_g(e,r)}$ is constructed via the
quotient:
\begin{equation}
\label{tpg}
Quot(U_H\rarr H, \bold{C}^n\otimes \oh_{U_{H}}) / SL_{N+1}\times
SL_n.
\end{equation}
There is a projection of the rigidified problem of
pairs $\{ (C,E)$ with isomorphisms (i) and (ii)$\}$ to a
rigidified moduli problem of curves $\{ C$ with isomorphism (i)$\}$.
The projection is $GL_{N+1}$-equivariant with respect to
the natural $GL_{N+1}$-action on the rigidified data of curves.
The Hilbert scheme $H$ is
a parameter space of the rigidified problem of curves (up to scalars).
By results of
Gieseker ([Gi]) reviewed in section (\ref{geese}),
the quotient $H/SL_{N+1}$ is
$\barr{M_g}$. A natural morphism $\barr{U_g(e,r)} \rarr \barr{M_g}$
is therefore obtained.
Gieseker's results require the choice in isomorphism (i) of at least
the $10$-canonical series.

The technical heart of the paper is the study of the Geometric
Invariant Theory problem (\ref{tpg}).
The method is divide and conquer. The action of
$SL_n$ alone is first studied. The $SL_n$-action is called the
fiberwise G.I.T. problem. It is solved in sections (\ref{puppy}) -
(\ref{lena}).
The action of $SL_{N+1}$ alone is then considered. There are
two pieces in the study of the $SL_{N+1}$ action. First, Gieseker's results
in [Gi] are used in an essential way. Second, the abstract
G.I.T. problem of
$SL_{N+1}$ acting on $\proj(Z)\times \proj(W)$ where $Z$, $W$
are representations of $SL_{N+1}$ is studied. If the linearization
is taken to be $\oh_{\proj(Z)}(k)\otimes \oh_{\proj(W)}(1)$ where
$k>>1$, there are elementary set theoretic
relationships between the stable and
unstable loci of $\proj(Z)$ and $\proj(Z)\times \proj(W)$. These
relationships are determined in section (\ref{abgit}).  Roughly speaking,
the abstract lemmas are used to import the invariants Gieseker
has determined in [Gi] to the problem at hand.
In section (\ref{conlo}), the
solution of fiberwise problem is combined with the study of the
$SL_{N+1}$-action
to solve the product G.I.T. problem (\ref{tpg}).

\subsection{Relationship with Past Results}
In [G-M],
D. Gieseker and I. Morrison propose
a different approach to a compactification over
 $\barr{M_g}$ of the universal
moduli space of slope-semistable bundles.
The moduli problem of pairs is rigidified by adding only the
data of an isomorphism
$\bold{C}^n \rarr H^0(C,E)$.
By further assumptions on $E$,
an embedding into a Grassmanian is obtained
 $$C \hookrightarrow \grass(rank(E),
H^0(C,E)^*).$$ The rigidified data is thus parametrized by a
Hilbert scheme of a Grassmanian. The $GL_n$-quotient problem is
studied to obtain a moduli space of pairs.

Recent progress along this alternate path has been made by
D. Abramovich, L. Caporaso,  and M. Teixidor ([A], [Ca], [T]).
A compactification, $\barr{P_{g,e}}$, of the universal
 Picard variety
is constructed in [Ca]. There is a natural isomorphism
$$\nu: \barr{P_{g,e}} \rarr \barr{U_g(e,1)}.$$
This isomorphism is established in section (\ref{lcap}).
In the rank $2$  case, the approach of [G-M] yields
a compactification not equivalent to $\barr{ U_g(e,2)}\ $ ([A]).
The higher rank constructions of [G-M] have not been completed.  They
are certainly expected to differ from $\barr{U_g(e,r)}$.

\subsection{Acknowledgements}
The results presented here constitute the author's 1994 Harvard
doctoral thesis.
It is a pleasure to
thank D. Abramovich and J. Harris for
introducing the author to
the higher rank compactification
problem. Conversations
with S. Mochizuki  on issues
both theoretical and technical have been of enormous
aid. The author has also benefited from discussions
with L. Caporaso, I. Morrison, H. Shahrouz, and M. Thaddeus.

\section{The Quotient Construction}
\subsection{Definitions}
\label{sfor}
Let $C$ be a genus $g\geq2$, Deligne-Mumford stable
curve. A coherent sheaf $E$ on $C$ is
 {\em torsion free} if $$\forall x\in C, \ \  depth_{\oh_x}(E_x)=1,$$ or
equivalently, if there does not exist a subsheaf
$$0\rarr F \rarr E $$
such that $dim(Supp(F))=0$. Let
$$C=\bigcup_{1}^{q}C_i $$
where the curves $C_i$ are the irreducible components of $C$.
Let $\om_i$ be the degree of the restriction of
the canonical bundle $\om_C$ to $C_i$.
Let
$r_i$ be the  the rank of $E$ at the generic point of $C_i$.
The {\em multirank} of $E$ is the $q$-tuple $(r_1, \ldots, r_q)$.
$E$ is of {\em uniform rank} $r$ if $r_i=r$ for each
$C_i$. If $E$ is of uniform rank $r$, define the degree of
$E$ by
$$e=\chi(E) - r(1-g).$$
A torsion free sheaf $E$ is defined
 to be {\em slope-stable (slope-semistable)}
if for each nonzero, proper subsheaf
$$0 \rarr F \rarr E $$
with multirank $(s_1, \ldots, s_q)$, the following inequality holds:
\begin{equation}
\label{bhab}
{\chi(F)\over \sum_{1}^{q}s_i \om_i} <
{\chi(E)\over \sum_{1}^{q}r_i\om_i},
\end{equation}
respectively,
$$\paren{{\chi(F)\over \sum_{1}^{q}s_i \om_i} \leq
{\chi(E)\over \sum_{1}^{q}r_i\om_i}}.$$
The above is Seshadri's definition of slope-(semi)stability
in the case of canonical polarization. In case $E$ is a vector
bundle on a nonsingular curve $C$, the slope-(semi)stability
condition (\ref{ahab}) of section (\ref{cmp}) and condition (\ref{bhab})
above coincide.
A slope-semistable sheaf has a
Jordan-Holder filtration with slope-stable factors.
Two
slope-semistable sheaves are {\em  equivalent} if
they possess the same set of Jordan-Holder factors.
Two equivalence classes are said to be {\em aut-equivalent} if
they differ by an automorphism of the underlying curve $C$.

For $g\geq2$ and  each pair of integers  $(e, \ r\geq1)$,
a projective variety $\barr{U_g(e,r)}$ and a
morphism
$$\eta: \barr{U_g(e,r)} \rarr \mgb$$
satisfying the following properties are constructed in
Theorem (\ref{priya}).
There is a functorial, bijective correspondence between the points
of  $\barr{U_g(e,r)}$ and aut-equivalence classes of
slope-semistable, torsion free
 sheaves of uniform rank $r$ and degree $e$ on
Deligne-Mumford stable curves of genus $g$.
The image of an aut-equivalence class under $\eta$ is the moduli
point of the underlying curve.

$\barr{U_g(e,r)}$ and $\eta$ will be
constructed via Geometric Invariant Theory.  The G.I.T. problem
is described in sections (\ref{geese}-\ref{linto}).
The solution is developed in sections (\ref{puppy}-\ref{conlo}) of the
paper.
Basic properties of
$\barr{U_g(e,r)}$ are studied in section (\ref{laslo}).
In particular, the equivalence of $\barr{U_g(e,1)}$ and
$\barr{P_{g,e}}$ is established in section (\ref{ender}).

\subsection{Gieseker's Construction}
\label{geese}
We review Gieseker's beautiful construction of $\barr{M_g}$.
Fix a genus $g \geq 2$. Define:
$$ d= 10(2g-2) $$
$$ N=d-g .$$
Consider the Hilbert scheme $H_{g,d,N}$
of genus $g$, degree $d$, curves in
$\proj^N_{\bold{C}}$. Let
$$H_g \subset H_{g,d,N}$$ denote the locus of nondegenerate,
$10$-canonical,
 Deligne-Mumford stable curves of genus g. $H_g$ is naturally a
closed subscheme of the open locus of nondegenerate, reduced, nodal
curves. In fact,
$H_g$ is
a nonsingular, irreducible, quasi-projective variety ([Gi]).
The symmetries of $\proj^N$ induce a natural $SL_{N+1}(\bold{C})$-action
on $H_g$.
D. Gieseker has studied the quotient $H_g/SL_{N+1}$ via
geometric invariant theory. It is shown in [Gi] that,
 for suitable linearizations,
$H_g/SL_{N+1}$ exists as a G.I.T. quotient and is
isomorphic to $\mgb$.

\subsection{Relative Quot Schemes}
\label{rqs}
Let $U_H$ be  the universal curve over
$H_g$.  We have a closed immersion
$$ U_H \hookrightarrow H_g \times \proj^N  $$ and
two projections:
$$\mu : U_H \rarr H_g$$
$$\nu: U_H \rarr \proj^N .$$
Let $\oh_U$ be the structure sheaf of $U_H$.
The
Grothendieck relative Quot scheme is central to our
construction.  We will be interested
in relative Quot schemes of the form
\begin{equation}
\label{qu}
Quot(\mu:U_H \rarr H_g,\  \tp \oh_U,\  \nu^*(\oh_{\proj^N}(1)),\  f)
\end{equation}
where $f$ is a Hilbert polynomial
with respect to the $\mu$-relatively very ample line
bundle $ \nu^*(\oh_{\proj^N}(1))$.
We denote the Quot scheme in (\ref{qu}) by
 $Q_g(\mu,n,f)$.

We recall the basic properties of the Quot scheme.
There is a canonical projective morphism
$ \pi: \Q \rarr H_g.$
The fibered product $\Q \times_{H_g} U_H$
is equipped with two projections:
$$ \theta : \Q \times_{H_g} U_H \rarr \Q$$
$$ \phi : \Q \times_{H_g} U_H \rarr U_H $$
and a  universal $\theta$-flat quotient
\begin{equation}
\label{uny}
\tp \oh_{Q\times U} \simeq \phi^* (\tp \oh_U) \rarr \E
\rarr 0.
\end{equation}
Let $\xi$ be a (closed) point of $\Q$.
The point $\pi(\xi) \in H_g$
corresponds to the $10$-canonical, Deligne-Mumford stable
curve $U_{\pi(\xi)}$. Restriction of the universal quotient
sequence (\ref{uny}) to $U_{\pi(\xi)}$ yields
a quotient
$$\tp \oh_{U_{\pi(\xi)}}\rarr \E_\xi \rarr 0$$ with
Hilbert polynomial
$$f(t)=\chi(\E_\xi\otimes \oh_{U_{\pi(\xi)}}(t)).$$
The above is a functorial bijective correspondence  between points
$\xi \in \Q$ and quotients
of $\tp \oh_{U_{\pi(\xi)}}$, $\pi(\xi) \in H_g$,  with Hilbert polynomial
$f$.

\subsection{Group Actions}
Denote the natural actions of $SL_{N+1}$ on $H_g$ and
$U_H$ by:
\begin{equation*}
\begin{CD}
U_H\times SL_{N+1} @>{a_U}>> U_H\\
@VV{\mu\times id}V @VV{\mu}V \\
H_g \times SL_{N+1} @>{a_H}>> H_g
\end{CD}
\end{equation*}
Also define:
$$\barr{\mu}: U_H \times SL_{N+1}
\stackrel{\mu \times inv}{\longrightarrow} H_g \times SL_{N+1}
\stackrel{a_H} {\longrightarrow} H_g.$$
$$\barr{\pi}:\Q \times SL_{N+1}
\stackrel{\pi \times id}{\longrightarrow} H_g \times SL_{N+1}
\stackrel{a_H}{\longrightarrow} H_g.$$
There is a natural isomorphism between the two
fibered products
$$\big(\Q \times SL_{N+1}\big) \times _{H_g} U_H \simeq
\Q \times _{H_g} \big( U_H \times SL_{N+1}\big)$$
where the projection maps to $H_g$ are $(\barr{\pi}, \mu)$ and
$(\pi, \barr{\mu})$ in the first and second products
respectively. The inversion in the definition of $\barr{\mu}$ is
required for the isomorphism of the fibered products.
There is a natural commutative diagram
\begin{equation*}
\begin{CD}
U_H\times SL_{N+1} @>{\barr{a}_U}>> U_H\\
@VV{\barr{\mu}}V @VV{\mu}V \\
H_g @>{id}>> H_g
\end{CD}
\end{equation*}
where $$\barr{a}_U: U_H \times SL_{N+1}
\stackrel{id \times inv} {\longrightarrow} U_H \times SL_{N+1}
\stackrel{a_U} {\rarr} U_H.$$
We therefore obtain a natural map of schemes over $\bold{C}$:
$$\varrho:(\Q \times SL_{N+1}) \times _{H_g} U_H \rarr \Q \times_{H_g} U_H.$$
By the functorial properties of $\Q$, the
$\varrho$-pull-back of the universal quotient sequence (\ref{uny}) on
$\Q \times_{H_g} U_H$ yields a natural group action:
$$\Q \times SL_{N+1} \rarr \Q .$$
Hence,
the natural $SL_{N+1}$-action
on $H_g$ lifts naturally to $Q_g(\mu,n,f)$.  There is
a natural $SL_n(\bold{C})$-action on $Q_g(\mu,n,f)$ induced by the
$SL_n(\bold{C})$-action on the tensor product $\tp \oh_U$.  In fact,
the $SL_{N+1}$-action
and the $SL_n$-action commute on $Q_g(\mu,n,f)$.
The commutation is most easily seen in the explicit linearized
projective embedding developed below in section (\ref{said}).
Hence, there  exists
a well-defined  $SL_{N+1} \times SL_n$-action.
For suitable choices of
$n$, $f$, and linearization, a component of the
quotient $Q_g(\mu,n,f)/(SL_{N+1} \times SL_n)$
will be $\barr{U_g(e,r)}$.

\subsection{Relative Embeddings}
\label{van}
Following [Gr], a family of
relative projective embeddings of $\Q$ over $H_g$ is constructed.
Since the inclusion
$$\Q \times_{H_g} U_g \hookrightarrow  \Q \times \proj ^N$$
is a closed immersion,
the universal quotient $\E$  can be extended by zero to
$\Q \times \proj^N$.
Let $$\tat: \Q \times\proj^N \rarr \Q$$ be the
projection.  The universal quotient sequence (\ref{uny})
induces the following sequence on $\Q \times \proj^N$:
$$ 0 \rarr \K \rarr \tp \oh_{Q \times \proj^N} \rarr \E \rarr 0 .$$
Since $\E$  and $\oh_{Q \times \proj^N}$ are $\tat$-flat,
$\K$ is $\tat$-flat. By the semicontinuity theorems for $\tat$-flat,
coherent sheaves, there exists an integer $t_{\al}$ such that for
each $t >t_{\al}$ and each $\xi \in \Q$ :
\begin{equation}
\label{cat}
h^1(\proj^N,\K_{\xi}\otimes \oh_{\proj^N}(t))=0
\end{equation}
\begin{equation}
\label{dog}
h^0(\proj^N,\E_{\xi}\otimes \oh_{\proj^N}(t))=f(t)
\end{equation}
\begin{equation}
\label{skunk}
h^1(\proj^N,\E_{\xi}\otimes \oh_{\proj^N}(t))=0
\end{equation}
\begin{equation}
\label{ppppp}
\tp H^0(\proj^N,\oh_{\proj^N}(t)) \rarr H^0(\proj^N, \E_{\xi}\otimes
\oh_{\proj^N}(t)) \rarr 0.
\end{equation}
The surjection of (\ref{ppppp})
follows from (\ref{cat}) and the
long exact cohomology sequence.
For each $t >t_{\al}$ there is a
well defined algebraic morphism (on points)
$$i_t: \Q \rarr \grass(f(t),
(\tp H^0(\proj^N, \oh_{\proj^N}(t)))^*)$$
defined by sending $\xi \in \Q$ to the subspace
$$ H^0(\proj^N, \E_\xi \otimes \oh_{\proj^N}(t))^* \subset
(\tp H^0(\proj^N,\oh_{\proj^N}(t)))^*.$$
By the theorems of Cohomology and Base Change,
it follows from conditions (\ref{cat}-\ref{skunk})
there exists
a surjection
$$\tp  H^0(\proj^N,\oh_{\proj^N}(t))\otimes
\oh_Q \simeq \theta_{\proj *}(\tp \oh_{Q \times \proj^N}(t))
\rarr \theta_{\proj *}(\E\otimes \oh_{\proj^N}(t)) \rarr 0$$
where
$\theta_{\proj *}(\E\otimes \oh_{\proj^N}(t))$ is a locally free, rank $f(t)$
quotient.
  The universal property of the Grassmanian
defines $i_t$ as a morphism of schemes.
It is known that there exists an integer $t_{\bay}$ such that
for all $t > t_{\bay}$, the morphism $\pi\times i_t$ is a closed
embedding:
$$\pi \times i_t: \Q \rarr H_g \times \grass(f(t),
(\tp H^0(\proj^N, \oh_{\proj^N}(t)))^*).$$
 The morphisms $\pi \times i_t, \ t>t_{\bay}(g,n,f)$,
form a countable family of relative projective embeddings
of the Quot scheme $Q_g(\mu,n,f)$ over $H_g$.

\subsection{Gieseker's Linearization}
\label{said}
Since the Hilbert scheme $H_{d,g,N}$ is the Quot scheme
$$Quot(\proj^N \rarr Spec(\bold C),\  \oh_{\proj^N}, \oh_{\proj^N}(1),\
 h(s)=ds-g+1),$$
there
are closed embeddings for $s >s_{\al}$:
$$i'_s: H_{d,g,N} \rarr \grass(h(s),
 H^0(\proj^N,\oh_{\proj^N}(s))^*).$$
By results of Gieseker, an integer $\barr{s}(g)$ can be chosen so
that the $SL_{N+1}$-linearized G.I.T. problem determined by $i'_{\barr{s}}$ has
two properties:
\begin{enumerate}
\item [(i.)] $H_g$ is contained in the stable locus.
\item [(ii.)] $H_g$ is closed in the semistable locus.
\end{enumerate}
In order to make use of (i) and (ii) above, we will only consider
immersions of the type $i'_{\barr{s}}$.

For each large $t$, there exists
an immersion:
$$i_{\barr{s},t}: \Q \rarr \grass(h(\barr{s}),
 H^0(\proj^N,\oh_{\proj^N}(\barr{s}))^*)
\times
\grass(f(t),
(\tp H^0(\proj^N,\oh_{\proj^N}(t)))^*). $$
By the Pl\"ucker embeddings, we obtain
$$j_{\barr{s},t}: \Q \rarr
 \proj(\bigwedge^{h(\barr{s})}H^0(\proj^N,\oh_{\proj^N}(\barr{s}))^*)
\times \proj(\bigwedge^{f(t)}(\tp H^0(\proj^N,\oh_{\proj^N}(t)))^*).$$
The fact that the $SL_{N+1}$ and $SL_{n}$-actions commute
on $\Q$ now follows from the observation that
the $SL_{N+1}$ and $SL_{n}$-actions commute on
$ \tp H^0(\proj^N,\oh_{\proj^N}(t))^*$.

\subsection{The G.I.T. Problem}
\label{linto}

Let $C$ be a Deligne-Mumford stable curve of genus $g\geq 2$.
For any coherent sheaf $F$ on $C$, it is not hard to see:
\begin{equation}
\label{oiler}
\chi(F\otimes \om_C^t) = \chi(F) + (\sum_{1}^{q}s_i\om_i)\cdot t
\end{equation}
where $(s_1,\ldots,s_q)$ is the multirank of $F$ ([Se]).
Equation (\ref{oiler}) and the slope inequalities
of section (\ref{sfor}) yield
a natural correspondence
$$(C,E) \rarr (C, E \otimes \om_C^t)$$
between slope-semistable, uniform rank $r$, torsion free sheaves
of degrees $e$ and $e+rt(2g-2)$. Therefore, it suffices to construct
$\barr{U_g(e,r)}$ for $e>>0$.

The strategy for studying the G.I.T. quotient
$$Q_g(\mu, n,f)/(SL_{N+1} \times SL_n)$$ is as follows.
First a rank $r\geq1$ is chosen.  Then the degree $e>e(g,r)$ is chosen
very large. The Hilbert polynomial
is determined by:
\begin{equation}
\label{goat}
f_{e,r}(t)= e+ r(1-g) + r10(2g-2)t.
\end{equation}
For $[C]\in H_g$,
$f_{e,r}(t)$ is
the Hilbert polynomial of degree $e$, uniform rank $r$, torsion free
sheaves on $C$ with respect to $\oh_{\proj^N}(1)$.
The integer $n$ is fixed by the Euler characteristic,  $n=f_{e,r}(0)$.
As remarked in section (\ref{moc})
of the introduction, $n$ will equal $h^0(C,E)$
for semistable pairs.
Let
$$\hat{t}(g,r,e)\ =\ t_{\bay}(g,n=f_{e,r}(0), f_{e,r})$$
be the
constant defined in section (\ref{van}) for $Q_g(\mu, f_{e,r}(0), f_{e,r})$.
A very large $t>\hat{t}(g,r,e)$
is then chosen.  Selecting $e$ and $t$ are the essential choices that
make the G.I.T. problem well-behaved.  Finally, to
determine a linearization of the
$SL_{N+1}\times SL_{n}$-action on the image of  $j_{\barr{s},t}$, weights
must be chosen on the two projective spaces.  These are chosen
so that almost all the weight is on the first,
\begin{equation}
\label{ffact}
\proj(\bigwedge^{h(\barr{s})}H^0(\proj^N,\oh_{\proj^N}(\barr{s}))^*).
\end{equation}
Since $SL_n$ acts only on the second factor of the product, the weighting
is irrelevant to the G.I.T. problem for the $SL_n$-action alone. The $SL_n$
action is studied in sections (\ref{puppy})-(\ref{lena}).
Since $SL_{N+1}$ acts on both factors,
the weighting is very relevant to the $SL_{N+1}$- G.I.T. problem. General
results of section (\ref{abgit})
show that in the case of extreme weighting, information
on the stable and unstable loci of the $SL_{N+1}$-action on
first factor can be transferred to the $SL_{N+1}$-action on the
product of the factors. Gieseker's study of the $SL_{N+1}$-action on the
first factor (\ref{ffact}) can therefore be used. In section (\ref{conlo}),
knowledge
of the $SL_n$ and $SL_{N+1}$ G.I.T problems is combined to solve the
$SL_{N+1}\times SL_n$ G.I.T.
problem on $\Q$.

\section{The Fiberwise G.I.T. Problem}
\label{puppy}
\subsection{The Fiberwise Result}
The fiber of $\pi :\Q \rarr H_g$ over a point $[C] \in H_g$ is
the Quot scheme
 $$Q_g(C,n,f)=Quot(C \rarr Spec(\bold{C}),\  \tp \oh_C ,\ \om_C^{10},\ f).$$
For large $t$, the morphism $i_t$ embeds $Q_g(C,n,f$) in
$$\grass(f(t), (\tp Sym^t(H^0(C,\om_C^{10})))^*).$$
The embedding $i_t$ yields an $SL_n$-linearized G.I.T. problem
on $Q_g(C,n,f).$ Before examining the global G.I.T. problem
for the construction of $\barr{U_g(e,r)},$ we will study this
fiberwise G.I.T. problem. The main result is:
\begin{tm}
\label{fred}
Let $g\geq2 $, $r>0$ be integers. There exist bounds $e(g,r)>r(g-1)$ and
$t(g,r,e)>\hat{t}(g,r,e)$ such that for each pair  $e>e(g,r)$,
$t>t(g,r,e)$ and any $[C]\in H_g$, the following holds:

\noindent
 A point $\ \xi \in Q_g(C,n=f_{e,r}(0),f_{e,r})$
corresponding to a
quotient $$ \tp \oh_C \rarr E \rarr 0$$
is G.I.T. stable (semistable) for the $SL_n$-linearization determined
by $i_t$  if and only if
 $E$ is a slope-stable (slope-semistable),  torsion free
sheaf on $C$ and
$$\psi:  \tp H^0(C,\oh_C) \rarr H^0(C,E)$$
is an isomorphism.
\end{tm}
P. Newstead has informed the author that
a generalization of this fiberwise G.I.T. problem
has been solved recently by C. Simpson in [Si]. A slight twist in our
approach is the uniformity of bound needed  for each
$[C] \in H_g$.  The proof will be developed  in many steps.

\subsection{The Numerical Criterion}
\label{nnmc}
Stability for a point $\xi$ in a linearized G.I.T. problem
can be checked by examining certain limits of $\xi$ along
$1$-parameter subgroups. This remarkable fact leads to the Numerical
Criterion. The most general form of the Numerical Criterion is presented in
section (\ref{nmc}). A  more precise version for the fiberwise G.I.T. problem
is stated here.

Fix a vector space $Z$ with the
trivial $SL_n$-action. In the applications below,
$$Z \stackrel{\sim}{=} Sym^t(H^0(C,\om_C^{10})).$$
Consider the linearized $SL_n$-action on $\grass(k, (\bold{C}^n \otimes
Z)^*)$
obtained from the standard representation of $SL_n$ on $\bold{C}^n$.
Let $\xi \in \grass(k,  (\bold{C}^n \otimes Z)^*).$
The element $\xi$ corresponds to a
$k$-dimensional quotient
$$ \rho_{\xi}: \bold{C}^n \otimes Z \rarr K_{\xi}.$$
Let $\barr{v} =( v_1, \ldots ,v_n)$ be a basis of
 $\bold{C}^n$ with integer
weights
$\barr{w} =(w(v_1), \ldots ,w(v_n) )$.
For combinatorial convenience,
the  additional condition that the weights sum to zero
is avoided here.
The representation weights of the corresponding $1$-parameter subgroup
of $SL_n$ are given by rescaling: $e_i=w(v_i)- \sum_{i} w(v_i)/n$.
An element $a \in \bold{C}^n \otimes Z$ is
said to be
{\em $\barr{v}$-pure} if it lies in a subspace of the form
$v_i\otimes Z$. The weight, $w(a)$,  of such an element is defined to be
$w(v_i)$.
The Numerical Criterion yields:
\begin{enumerate}
\item $\xi$ is unstable if and only if there exists a
basis $\barr{v}$ of $\bold{C}^n$ and weights
$\barr{w}$ with the
following property.  For any $k$-tuple of  $\barr{v}$-pure
elements $(a_1, \ldots, a_k ) $ such that
$(\rho_{\xi}(a_1), \ldots, \rho_{\xi}(a_k) )$ is
a basis of $K_{\xi}$, the inequality
$$\sum_{i=1}^{n} {w(v_i)\over n}  < \sum_{j=1}^{k}{w(a_j)\over k} $$
is satisfied.
\item$\xi$ is stable (semistable) if and only if
for every basis $\barr{v}$ of $\bold{C}^n$ and any nonconstant weights
$\barr{w}$ the following holds. There exist
$\barr{v}$-pure
elements $(a_1, \ldots, a_k )$ such that
$( \rho_{\xi}(a_1), \ldots, \rho_{\xi}(a_k) )$ is
a basis of $K_{\xi}$ and
$$\sum_{i=1}^{n} {w(v_i)\over n} \  >(\geq) \  \sum_{j=1}^{k}{w(a_j)\over k}.$$
\end{enumerate}
See, for example, [N] or [M-F].

\subsection{Step I}
The instability arguments will use the following Lemma.
\begin{lm}
\label{archy}
Let $g\geq2$, $r>0$, $e>r(g-1)$ be integers, $[C]\in H_g$.
Suppose $\xi \in Q_g(C,n=f_{e,r}(0),f_{e,r})$ corresponds to a
quotient $$ \tp \oh_C \rarr E \rarr 0.$$
Let $U\subset \bold{C}^n$ be a subspace. Let
$\psi(U\otimes H^0(C,\oh_C))= W \subset H^0(C,E)$.
Let $G$ be the subsheaf of $E$ generated by $W$.
For any $t>\hat{t}(g,e,r)$ the following holds: if
\begin{equation}
\label{zara}
{dim(U) \over n} > {h^0(C, G\otimes \om_C^{10t}) \over f_{e,r}(t)},
\end{equation}
then $\xi$ is G.I.T. unstable for the $SL_n$-linearization
determined by $i_t$.
\end{lm}
\begin{pf}
Let $u=dim(U)$.
Inequality (\ref{zara}) implies $0<u<n$. Let $\barr{v}$ be a basis
of $\bold{C}^n$ such that $(v_1,\ldots, v_u)$ is a basis of $u$.
Select weights as follows: $w(v_i)=0$ for $1\leq i \leq u$ and
$w(v_i)=1$ for $u+1\leq i \leq n$.  We now use the Numerical Criterion
for the $SL_n$-action on
$\grass(f_{e,r}(t),(\tp Sym^t(H^0(C,\om_C^{10})))^*)$. The element
$\xi$ corresponds to a quotient
$$\psi^t:\tp Sym^t(H^0(C,\om_C^{10})) \rarr H^0(C,E\otimes \om_C^{10t})
\rarr 0.$$
Suppose $(a_1,\ldots, a_{f_{e,r}(t)})$ is a tuple of $\barr{v}$-pure
elements mapped by $\psi^t$ to a basis of $H^0(C,E\otimes \om_C^{10t})$.
All $a_j$'s have weight $1$ except those contained in $U\otimes
Sym^t(H^0(C,\om_C^{10}))$ which have weight $0$. The number
of $a_j$'s of weight $0$ is hence bounded by $h^0(C,G\otimes \om_C^{10t})$.
Since
$$\sum_{i=1}^{n}{w(v_i)\over n}= 1-{u\over n} < 1-
{h^0(C,G\otimes \om_C^{10t}) \over f_{e,r}(t)} \leq \sum_{j=1}^{f_{e,r}(t)}
{w(a_j)\over f_{e,r}(t)},$$
the Numerical Criterion implies $\xi$ is unstable.
\end{pf}

\begin{pr}
\label{wk1}
Let $g\geq2$, $r>0$ be integers. For each pair $e>r(g-1)$,
$t>\hat{t}(g,r,e)$ and any $[C]\in H_g$, the following holds:

\noindent
If $\xi\in Q_g(C,n=f_{e,r}(0) f_{e,r})$ corresponds to a quotient
$$\tp \oh_C \rarr E \rarr 0$$
where $\psi: \tp H^0(C,\oh_C) \rarr H^0(C,E)$ is not injective,
then $\xi$ is G.I.T. unstable for the $SL_n$-linearization determined
by $i_t$.
\end{pr}

\begin{pf}
Suppose $\psi$ is not injective.
Let $U\otimes H^0(C,\oh_C)$ be the nontrivial kernel of $\psi$.
The assumptions of Lemma (\ref{archy}) are easily checked since
$W=0$ and $G$ is the zero sheaf. $\xi$ is G.I.T unstable by
Lemma (\ref{archy}).
\end{pf}

\section{Cohomology Bounds}
\label{bondy}
\subsection{The Bounds}
In order to further investigate the fiberwise $SL_n$-action,
we need to control the first cohomology  in various ways.
\begin{lm}
\label{bat}
Let $g\geq 2$, $R>0$ be integers. There exists an
integer $b(g,R)$ with the following property. If
 $E$ is a coherent sheaf
on a Deligne-Mumford stable, genus $g$ curve $C$ such that:
\begin{enumerate}
\item[(i.)] $E$ is generated
by global sections.
\item[(ii.)] $E$ has generic rank less than $R$ on each irreducible
component of $C$.
\end{enumerate}
Then
$h^1(C,E) < b(g,R).$
\end{lm}

\begin{pf}
Since $\om_C$ is ample and of degree $2g-2$,
there is a bound $q(g)=2g-2$ on the number of components of $C$.
Since $E$ is generated by global sections and has bounded
rank, there exists an exact  sequence:
$$ \sumo_{1}^{qR}\oh_C \rarr E \rarr \tau \rarr 0 $$
where $Supp(\tau)$ is at most dimension zero.
Hence
$$ h^1(C,E) \leq qR \cdot h^1(C,\oh_C).$$
Since $h^1(C,\oh_C)=g$,
$\ b(g,R)= (2g-2)Rg +1$ will have the
required property.
\end{pf}

\begin{lm}
\label{cow}
Let $g\geq 2$, $R>0$, be integers.  Let $E$ be a coherent
sheaf on a Deligne-Mumford stable, genus $g$ curve $C$
satisfying (i) and (ii) of Lemma (\ref{bat}).
Suppose $F$
is a subsheaf of $E$ generated by global sections.
Then
 $$|\chi(F)| <  | \chi(E)| +b(g,R) .$$
\end{lm}

\begin{pf}
By Lemma (\ref{bat}),
$h^1(C,F)<b(g,R)$.  Therefore,  $- b(g,R) < \chi(F)$. By
Lemma (\ref{bat}) applied to $E$,
$$\chi(F) \leq h^0(C,F) \leq h^0(C,E) < \chi(E) + b(g,r) \leq
|\chi(E)| + b(g,r).$$ The result follows.
\end{pf}

\begin{lm}
\label{cotton}
Let $g\geq 2$, $R>0$, $\chi$ be integers.
There exists an
integer $p(g,R,\chi)$ with the following
property. Let $E$ be any coherent
sheaf on any Deligne-Mumford stable, genus $g$ curve $C$
satisfying (i) and (ii) of Lemma (\ref{bat}) and satisfying
$\chi(E)=\chi$.
Let $F$
be any subsheaf of $E$ generated by $k$ global sections:
\begin{equation}
\label{eee}
\bold{C}^k \otimes \oh_C \rarr F \rarr 0.
\end{equation}
Then
for all $t > p(g,R,\chi)$ :
\begin{enumerate}
\item[(i.)] $h^1(C,F\otimes \om_C^{10t})=0.$
\item[(ii.)] $\bold{C}^k \otimes Sym^t(H^0(C,\om_C^{10}))
 \rarr H^0(C,F\otimes \om_C^{10t})  \rarr 0.$
\end{enumerate}
\end{lm}

\begin{pf}
Let $C=\bigcup_{1}^{q}C_i$. Let
$\om_i$ be the degree of $\om_C$ restricted to $C_i$.
Let $(s_1, \ldots, s_q)$ be the multirank
of $F$.  By (\ref{oiler}) of section (\ref{linto}), the Hilbert polynomial
of $F$ with respect to $\om_C^{10}$ is:
$$\chi(F\otimes \om_C^{10t})= \chi(F)+ (\sum_{1}^{q}s_i \om_i)
\cdot 10t. $$
By Lemma (\ref{cow}), $|\chi(F)| < |\chi| +b(g,R)$.
Also
$$0\leq s_i <R,\ \ \  1\leq q \leq 2g-2,\ \ \  1\leq \om_i \leq 2g-2.$$
Therefore the data $g$, $R$, and $\chi$
determine a finite collection of Hilbert polynomials
$$\{f_1, \ldots , f_m \}$$
 that contains
the Hilbert polynomial of every allowed sheaf $F$.

The morphism (\ref{eee}) yields a
 natural map $\psi: \bold{C}^k \rarr H^0(C,F).$
Let $$im(\psi) = V\subset H^0(C,F).$$
We note that (\ref{eee}) can be factored:
$$\bold{C}^k \otimes \oh_C  \rarr V\otimes \oh_C \rarr F \rarr 0.$$
Since (ii) is surjective if and only if the analogous map in which
$\bold{C}^k$ is replaced by $V$ is surjective,
we can assume
$$k\leq h^0(C,F) \leq  h^0(C,E) < |\chi| + b(g,R).$$

Suppose $F$
is a coherent sheaf on a Deligne-Mumford
stable, genus $g$ curve satisfying:
\begin{enumerate}
\item [(a.)] $F$ is generated by $k< |\chi|+ b(g,R)$
 global sections: $\bold{C}^{k}\otimes \oh_C \rarr F \rarr 0.$
\item [(b.)] $F$ has Hilbert polynomial
$f$ (with respect to $\om_C^{10}$).
\end{enumerate}
Then there exists a integer ${\barr{t}}(g,k,f)$ such that
for all $t > {\barr{t}}(g,k,f)$ :
\begin{enumerate}
\item[(i.)] $h^1(C,F\otimes \om_C^{10t})=0.$
\item[(ii.)] $\bold{C}^{k}\otimes  Sym^t(H^0(C,\om_C^{10}))
 \rarr H^0(C,F\otimes \om_C^{10t})  \rarr 0.$
\end{enumerate}
The
existence of ${\barr{t}}(g,k,f)$ follows from
statements  (\ref{skunk}) and (\ref{ppppp}) of section
(\ref{van}) applied to the Quot scheme $Q_g(\mu,k,f)$.
Now let
 $$p(g,R,\chi)=max \{\ {\barr{t}}(g,k,f_j)
 \ | \ 1\leq k \leq |\chi|+b(g,R), \ 1\leq j \leq m \} .$$
It follows easily that $p(g,R,\chi)$ has the required property.
\end{pf}

\subsection{Step II} We apply these cohomology bounds along with the
Numerical Criterion in another simple case. First, two definitions:
\label{rdef}

Define  $R(g,r)=r(2g-2)+1$.
If $E$ is a coherent sheaf on $C$ with  multirank $(r_i)$ and
Hilbert polynomial $f_{e,r}$
with respect to $\om_C^{10}$, then by
(\ref{oiler}) of section (\ref{linto}), $$\sum r_i \om_i = r(2g-2).$$
Therefore, $r_i < R(g,r)$ for each $i$.

If $E$ is a coherent sheaf on $C$, there is canonical
sequence
$$0 \rarr \tau_E \rarr E \rarr E' \rarr 0$$
where $\tau_E$ is the {\em torsion subsheaf} of $E$  and $E'$ is torsion free.

\begin{pr}
\label{bob}
Let $g\geq2 $, $r>0$, $e>r(g-1)$ be integers. There exists a bound
$t_0(g,r,e)>\hat{t}(g,r,e)$ such that for each
$t > t_0(g,r,e)$, and $[C]\in H_g$ the following holds:

\noindent
If $\ \xi \in Q_g(C,n=f_{e,r}(0),f_{e,r})$ corresponds to
a quotient
\begin{equation}
\label{deer}
\tp \oh_C \rarr E \rarr 0
\end{equation}
where
$\psi\bigl( \tp H^0(C,\oh_C) \bigr) \cap H^0(C,\tau_E) \neq 0$,
then $\xi$ is G.I.T. unstable for the $SL_n$-linearization
determined by $i_t$.
\end{pr}

\begin{pf}
Let $U\subset \bold{C}^n$ be a $1$ dimensional subspace such that
$$\psi(U\otimes H^0(C,\oh_C))=W \subset H^0(C,\tau_E).$$ Let $G$ be
the subsheaf of $E$ generated by $W$.
For all $t$,
$$h^0(C, G\otimes \om_C^{10t}) \leq
h^0(C,\tau_E \otimes \om_C^{10t})=h^0(C,\tau_E).$$
By Lemma (\ref{bat}),$$h^0(C,\tau_E) \leq h^0(C,E) =\chi(E)+ h^1(C,E) <
 f_{e,r}(0)+ b(g,R(g,r)).$$
There certainly
exists a $t_0(g,r,e) > \hat{t}(g,r,e)$ satisfying $\forall t > t_0(g,r,e)$,
$${1\over n} >
{f_{e,r}(0)+b(g,R(g,r)) \over f_{e,r}(t)}.$$
The Proposition is now a consequence of Lemma (\ref{archy}).
\end{pf}

\section{Slope-Unstable, Torsion Free Sheaves}
\label{ggh}
\subsection {Step III}
Propositions (\ref{wk1}) and (\ref{bob}) conclude G.I.T. instability from
certain undesirable properties of points in
$Q_g(C,n,f_{e,r})$. In this section,  G.I.T. instability is
concluded from slope-instability in the case
where $\psi$ is an isomorphism and $E$ is torsion free.
The case where $\psi$ is not an isomorphism (and $E$ is arbitrary) is
analyzed in section (\ref{nassy}) where G.I.T.
instability is established.
The above results (for suitable choices of constants and linearizations) show
only points of $Q_g(C,n,f_{e,r})$ where $\psi$ is an isomorphism and
$E$ is a slope-semistable torsion free sheaf may be G.I.T. semistable.
The G.I.T.
(semi)stability results are established in section (\ref{lena}).
\begin{pr}
\label{kiwi}
Let $g\geq 2$, $r>0$ be integers. There exist bounds
$e_1(g,r)>r(g-1)$ and $t_1(g,r,e)>\hat{t}(g,r,e)$ such
that
for each pair  $e>e_1(g,r)$, $t>t_1(g,r,e)$ and any
$[C]\in H_g$, the following holds:

\noindent
If $\ \xi \in Q_g(C,n=f_{e,r}(0),f_{e,r})$ corresponds to
a quotient
$$\tp \oh_C \rarr E \rarr 0$$
where
$\psi: \tp H^0(C,\oh_C) \rarr H^0(C,E)$
is an isomorphism and
$E$ is a slope-unstable, torsion free  sheaf,
then $\xi$ is G.I.T. unstable for the $SL_n$-linearization
determined by $i_t$.
\end{pr}

\subsection{Lemmas and Proof} The proof of Proposition (\ref{kiwi})
requires two Lemmas which are used to apply Lemma (\ref{archy}).
First a destabilizing subsheaf $F$ of $E$ is selected. $F$ determines
a filtration: $ W=H^0(C,F) \subset H^0(C,E)$. $H^0(C,E)$ is identified with
$\bold{C}^n$ by $\psi$. Let $U=\psi^{-1}(W)$.
If $F$ is generated by global sections,
the vanishing
theorems of section (\ref{bondy}) can be applied. Riemann-Roch then
shows the conditions of Lemma (\ref{archy})
for $U$, $W$ follow from the destabilizing
property of $F$ (Lemma (\ref{book})).
In fact, the vanishing argument is valid when
$F$ is generically generated by global sections. Lemma (\ref{water}) guarantees
that a destabilizing subsheaf $F$ generically generated by global sections
exists if $E$ is of high degree.

\begin{lm}
\label{water}
Let $g\geq2$, $r>0$ be integers.
There exists an integer $e_1(g,r)>r(g-1)$
such that for each
$e> e_1(g,r)$ and $[C]\in H_g$ the following holds:

\noindent
If $E$ be a slope-unstable, torsion free sheaf on
$C$ with Hilbert polynomial
$f_{e,r}$ (with respect to $\om_C^{10}$), then there
exists a nonzero, proper destabilizing subsheaf $F$ of $E$ and an
exact sequence
\begin{equation}
\label{monk}
0 \rarr \barr{F} \rarr F \rarr \tau \rarr 0
\end{equation}
where $\barr{F}$ is generated by global sections
and $Supp(\tau)$ is at most  dimension zero.
\end{lm}

\begin{pf}
Since $E$ is slope-unstable, there exists a nonzero, proper
destabilizing subsheaf,
$$0\rarr F \rarr E.$$
Let $C$ be the union of components $\{C_i \}$ where $1\leq i \leq q$.
and let $(s_i)$, $(r_i)$ be the multiranks
of $F$, $E$. Since $E$ is torsion free and $F$ is nonzero, the
multirank of $F$ is not identically  zero.
Since the Hilbert polynomial of $E$ is
$f_{e,r}$, we see (by section (\ref{rdef})) $R(g,r)= r(2g-2)+1$
satisfies
$\forall i \ \ r_i < R(g,r)$.
 $F$ can be chosen to have minimal multirank
in the following sense. If $F'$ is a nonzero
subsheaf of $F$ with
multirank $(s'_i )$ such that
$\exists j \ \  s'_j< s_j$ , then $F'$ is not destabilizing.
Let $\barr{F}$ be the subsheaf of $F$ generated by  the
global sections $H^0(C,F)$.
Since $F$ is destabilizing:
$$ h^0(C,F) \geq \chi(F) > \chi(E)\cdot \paren{\sum s_i \om_i\over r(2g-2)} =
(e+r(1-g)) \cdot \paren{\sum s_i \om_i\over r(2g-2)}.$$
Hence if $e >r(g-1)$, $h^0(C,F)>0$ and $\barr{F}$ is nonzero.
We now assume $e > r(g-1)$.
Let $(\barr{s}_i )$ be the nontrivial
multirank of $\barr{F}$.  The sequence
(\ref{monk}) has the required properties  if and only if
$(\barr{s}_i )=(s_i )$. Suppose $\exists j$,
$\ \barr{s}_j < s_j$. Then $\barr{F}$ is not destabilizing, so
$$\chi(\barr{F}) \leq
(e+r(1-g)) \cdot \paren{\sum \barr{s}_i \om_i\over r(2g-2)}.$$
We obtain
$$\chi(\barr{F}) <  h^0(C,F) \cdot \paren{\sum \barr{s}_i \om_i \over
 \sum s_i \om_i } \leq
h^0(C,F) \cdot  \paren{r(2g-2)-1 \over r(2g-2)}.$$
The last inequality follows from the fact
$$0<\sum \barr{s}_i \om_i < \sum s_i \om_i \leq r(2g-2).$$
Since $\barr{F}$ is generated by global sections, Lemma (\ref{bat})
yields
$$h^1(C,\barr{F}) < b(g,R(g,r))=b.$$
We conclude,
$$h^0(C,\barr{F}) <  b+ h^0(C,F) \cdot \paren{r(2g-2)-1\over r(2g-2)}.$$
Since $h^0(C,\barr{F}) =h^0(C,F)$ and
$$h^0(C,F) > {e+r(1-g)\over r(2g-2)}\ ,$$
we obtain the bound
$$\paren{e+r(1-g) \over r(2g-2)} \cdot \paren{ 1\over r(2g-2)}<
b.$$
Hence
$$e_1(g,r) = b(g,R(g,r)) \cdot (r^2(2g-2)^2)+r(g-1)$$
has the property required by the Lemma.
\end{pf}

\begin{lm}
\label{book}
Let $g \geq2$, $r>0$, $e>e_1(g,r)$ be integers. There
exists an integer $t_1(g,r,e)> \hat{t}(g,r,e)$ such that
for each
$t>t_1(g,r,e)$ and $[C]\in H_g$, the following holds:

\noindent
If $\xi \in Q_g(C, n=f_{e,r}(0),f_{e,r})$ corresponds to a
quotient
$$ \tp \oh_C \rarr E \rarr 0$$
where $E$ is a slope-unstable, torsion free sheaf on $C$,
 then there exists a nonzero, proper subsheaf
$$0\rarr F \rarr E$$
such that
\begin{equation}
\label{crow}
 {h^0(C,F)\over n}>
 {h^0(C,F\otimes \om_C^{10t})\over f_{e,r}(t)}\ .
\end{equation}
\end{lm}

\begin{pf} Let $t_1(g,r,e)>p(g,R(g,r),\chi=f_{e,r}(0))$ be
determined by Lemma (\ref{cotton}).
Take $F$ to be a nonzero, proper, destabilizing subsheaf of
 $E$ for which there
exists a  sequence
\begin{equation}
\label{god}
0 \rarr \barr{F} \rarr F \rarr \tau \rarr 0
\end{equation}
where $\barr{F}$ is generated by global sections and
$Supp(\tau)$ has dimension zero. Such $F$ exist by Lemma (\ref{water}). Since
$\barr{F}$ is a subsheaf of $E$ and is generated by global
sections, Lemma (\ref{cotton}) yields for any $t >t_1$,
$ h^1(C,\barr{F}\otimes \om_C^{10t})=0$.
By the exact sequence in cohomology of (\ref{god}),
$ h^1(C,F\otimes \om_C^{10t})=0$.
Let $(s_i)$ be the (nontrivial) multirank of $F$.
We have
$$h^0(C,F\otimes \om_C^{10t})= \chi(F) + (\sum s_i \om_i)10t,$$
$$f_{e,r}(t)= \chi(E) + r(2g-2)10t.$$
We obtain  $$\chi(F) \cdot f_{e,r}(t)-
 \chi(E) \cdot h^0(C,F\otimes \om_C^{10t}) = $$
$$\chi(F) \cdot r(2g-2)10t - \chi(E) \cdot (\sum s_i \om_i)10t >0.$$
The last inequality follows from the destabilizing property of $F$.
Hence
$$ {\chi(F)\over \chi(E)} >
 {h^0(C,F\otimes \om_C^{10t})\over f_{e,r}(t)}.$$
Since $h^0(C,F) \geq \chi(F)$ and $\chi(E)=n$, the
Lemma is proven.
\end{pf}

We can now apply Lemma (\ref{archy}).
\begin{pf}[Of Proposition (\ref{kiwi})]
Let $e_1(g,r)$ be as in Lemma (\ref{water}). For $e>e_1(g,r)$,
let $t_1(g,r,e)$ be determined by Lemma (\ref{book}). Suppose
$t>t_1(g,r,e)$.
Let
$F$ be the subsheaf of $E$ determined by Lemma (\ref{book}).
Let $U\subset \bold{C}^n = \psi^{-1} (H^0(C,F))$.
Since $\psi$ is an isomorphism $dim(U)=h^0(C,F)$. Let $G$ be the
subsheaf generated by the global sections $H^0(C,F)$. Certainly
$$h^0(C, F\otimes \om_C^{10t}) > h^0(C, G\otimes \om_C^{10t}).$$
Lemmas (\ref{book}) and (\ref{archy}) are now sufficient to conclude
the desired G.I.T. instability.
\end{pf}

\section {Special, Torsion Bounded Sheaves}
\label{nassy}
\subsection{Lemmas}
As always, suppose $n=f_{e,r}(0)= \chi (E)$.
If $$\psi: \bold{C}^{n}\otimes  H^0(C,\oh_C) \rarr H^0(C,E) $$
is injective but not
an isomorphism, then $h^1(C,E)\neq 0$.
We now investigate this case and conclude G.I.T.
instability for the corresponding points of
$Q_g(C,n,f_{e,r})$. The strategy is the following. Since
$H^1(C,E)$ is dual to $Hom(E,\om_C)$, the latter must be nonzero.
In Lemma (\ref{boy}),
The kernel of a nonzero element of $Hom(E, \om_C)$ is analyzed
to produce a very destabilizing subsheaf $F$ of $E$. Lemma (\ref{archy})
is then applied as in section (\ref{ggh}).
In order to carry out the above plan,
the torsion of $E$ must be treated with care.

For any
coherent sheaf $E$ on $C$, let
$0 \rarr \tau_E \rarr E $
be the torsion subsheaf.  $E$ is said to have {\em torsion bounded}
by $k$ if $\chi(\tau_E) < k.$
Let $R(g,r)=r(2g-2)+1$ as defined in section (\ref{rdef}).

\begin{lm}
\label{boy}
Let $g\geq2$, $r>0$ be integers.  There exists an integer
 $e_2(g,r)> r(g-1)$ such that for each
$e >e_2(g,r)$ and $[C]\in H_g$, the
following holds:

\noindent
If $E$ is  a  coherent sheaf
on $C$ with
Hilbert polynomial $f_{e,r}$ (with respect to $\om_C^{10}$) satisfying
\begin{enumerate}
\item[(i.)] $h^1(C,E) \neq 0$
\item[(ii.)] $E$ has torsion bounded by $b(g,R(g,r))$,
\end{enumerate}
then there
exists a nonzero, proper  subsheaf $F$ of $E$
with multirank $(s_i)$ not identically zero such that
\begin{enumerate}
\item [(i.)]$F$ is generated by global sections
\item [(ii.)]$${\chi(F) - b(g,R(g,r)) \over \sum s_i\om_i}>
{\chi(E)\over r(2g-2)} +1\ .$$
\end{enumerate}

\end{lm}

\begin{pf}
Since by Serre duality $H^1(C,E)^* \cong Hom(E,\om_C)$,
there exists a nonzero
morphism of coherent sheaves:
$$\sigma : E \rarr \om_C .$$
We have
$0 \rarr \sigma(E) \rarr \om_C$
where $\sigma(E) \neq 0$.
Since $\om_C$ is torsion free, $\sigma(E)$ has multirank
not identically zero. Note
$$\chi(\sigma(E)) \leq h^0(C,\sigma(E)) \leq h^0(C,w_C)=g .$$
Consider the exact sequence:
$$0 \rarr K \rarr E \rarr \sigma(E) \rarr 0.$$
Since $\chi(K)= \chi(E) - \chi(\sigma(E))$,
$$\chi(K) \geq \chi(E) - g .$$
For $e>r(g-1) +g $, $\chi(K) >0$ and $K \neq 0$.
Let $F$ be the subsheaf generated by the global sections of $K$.
$\chi(K) >0$ implies $F \neq 0$. Let $b=b(g,R(g,r))$.  We have
$$\chi(F)> h^0(C,F) - b = h^0(C,K)-b \geq \chi(K) -b
\geq \chi(E)-b-g .$$
For $e> r(g-1)+2b+g$, $\chi(F) > b$.
Now assume $e> r(g-1)+2b+g.$
By the bound on the torsion of $E$, $F$ is not contained
in $\tau_E$.
Let $(s_i)$ be the multirank of $F$. Since $F$ is not torsion,
the multirank is not identically zero. In fact,
since $\sigma(E)$ has multirank
not identically zero,
$$0< \sum s_i\om_i < r(2g-2).$$
We conclude
\begin{eqnarray*}
{\chi(F)-b \over \sum s_i \om_i} &>&
\paren{\chi(E)-2b-g \over r(2g-2)} \cdot \paren{r(2g-2)\over \sum s_i\om_i} \\
& \geq &
\paren{\chi(E)-2b-g \over r(2g-2)} \cdot \paren{ r(2g-2)\over r(2g-2)-1} \ .
\end{eqnarray*}
For large $e$ depending only on $g$ and $r$,
$$\paren{\chi(E)-2b-g \over r(2g-2)} \cdot \paren{r(2g-2) \over r(2g-2)-1}
> {\chi(E) \over r(2g-2)} +1\ .$$
We omit the explicit bound.
\end{pf}
An analogue of Lemma (\ref{book}) is now proven.
\begin{lm}
\label{bag}
Let $g \geq2$, $r>0$, $e>e_2(g,r)$, be integers. There
exists an integer $t_2(g,r,e)> t_0(g,r,e)$ such
for each
 $t > t_2(g,r,e)$ and $[C]\in H_g$, the following holds:

\noindent
If $\xi \in Q_g(C, n=f_{e,r}(0),f_{e,r})$
corresponds to a
quotient
$$ \tp \oh_C \rarr E \rarr 0$$
where $E$ is a coherent sheaf  on $C$ satisfying
\begin{enumerate}
\item[(i.)] $ \psi: \tp H^0(C,\oh_C) \rarr H^0(C,E) $
is injective
\item[(ii.)] $h^1(C,E)\neq 0$
\item[(iii.)] $E$ has torsion bounded by $b(g,(R,g,r))$,
\end{enumerate}
 then there exist a nonzero subspace $W \subset
\psi \Bigl( \tp H^0(C,\oh_C) \Bigr)$
generating a nonzero, proper subsheaf
$0\rarr G \rarr E$
such that
\begin{equation}
\label{snow}
{dim(W)\over n} >
 {h^0(C,G\otimes \om_C^{10t})\over f_{e,r}(t)}\  .
\end{equation}
\end{lm}

\begin{pf}
Let $F$ be the subsheaf of $E$ determined by Lemma (\ref{boy}).
Let $$W=im(\psi) \cap H^0(C,F) .$$
Since $$h^0(C,E)< \chi(E) + b(g,R(g,r)),$$
and $\psi$ is injective,
$$dim(W) > h^0(C,F) -b \geq \chi(F) -b.$$
Note by condition (ii) of $F$ in Lemma (\ref{boy}), $dim(W) > 0$.
Let $$t > p(g, R(g,r), \chi=f_{e,r}(0)).$$ Since
$F$ is generated by global sections,
$$h^1(C,F\otimes w_C^{10t})=0$$
by Lemma (\ref{cotton}).
We have
$$h^0(C,F\otimes \om_C^{10t}) =\chi(F) + (\sum s_i\om_i)10t,$$
$$f_{e,r}(t)=\chi(E) + r(2g-2)10t.$$
We compute
$$(\chi(F)-b)\cdot f_{e,r}(t)-
\chi(E) \cdot h^0(C,F\otimes \om_C^{10t}) =  $$
$$ (\chi(F)-b) \cdot r(2g-2)10t
 -\chi(E) \cdot (\sum s_i\om_i)10t -b\cdot\chi(E) >$$
$$ r(2g-2)\cdot (\sum s_i\om_i) \cdot 10t - b \cdot \chi(E).$$
The last inequality follows from condition (ii) of $F$ in
Lemma (\ref{boy}).
If also $$t>  b \cdot \chi(E)  = b(g,R(g,r)) \cdot (e+r(1-g)),$$
then
$$ {\chi(F)-b \over \chi(E)} >
 {h^0(C,F\otimes \om_C^{10t})\over f_{e,r}(t)}\  .$$
Let $G$ be the subsheaf of $F$ generated by $W$.
Since $dim(W) > \chi(F) -b $, $n=\chi(E)$, and
$$h^0(C,F\otimes \om_C^{10t}) \geq h^0(C,G\otimes \om_C^{10t}),$$
 the proof is complete.
\end{pf}

\subsection{Step IV}

\begin{pr}
\label{cake}
Let $g\geq 2$, $r>0$ be integers. There exist bounds
$e_2(g,r)>r(g-1)$ and $t_2(g,r,e)>t_0(g,r,e)$ such that
for each pair $e>e_2(g,r)$, $t>t_2(g,r,e)$ and any $[C]\in H_g$,
 the following holds:

\noindent
If $\ \xi \in Q_g(C,n=f_{e,r}(0),f_{e,r})$ corresponds to
a quotient
$$\tp \oh_C \rarr E \rarr 0$$
where
$h^1(E,C) \neq 0$,
then $\xi$ is G.I.T. unstable for the $SL_n$-linearization
determined by $i_t$.
\end{pr}

\begin{pf}
Let $e_2(g,r)$ be given by Lemma (\ref{boy}).
For $e>e_2(g,r)$, let $t_2(g,r,e)$ be given by
Lemma (\ref{bag}).
Let $$\psi : \tp H^0(C,\oh_C) \rarr H^0(C,E)$$
be the map on global sections.
If $\psi$ has a nontrivial kernel, $\xi$ is unstable by
Proposition (\ref{wk1}).  We can assume
$\psi$ is injective.
Note that $im(\psi)$ has codimension less than $b(g,R(g,r))$ in
$H^0(C,E)$. If
$0 \rarr \tau \rarr E $
is a torsion subsheaf such that $h^0(C,\tau)=\chi(\tau) \geq  b(g,R(g,r))$,
then $$im(\psi) \cap H^0(\tau,C) \neq 0.$$
In this case, since $t>t_0(g,r,e)$,
 $\xi$ is unstable by Proposition (\ref{bob}).
We can assume $E$ has torsion bounded by $b$.
We now can apply Lemma (\ref{bag}). Let
$W\subset im(\psi)$ be determined by Lemma (\ref{bag}).
Let $U=\psi^{-1}(W)$.
Since $\psi$ is injective, $dim(U)=dim(W)$.
Lemmas (\ref{bag}) and (\ref{archy}) now imply the desired
G.I.T. instability.
\end{pf}

\section{Slope-Semistable, Torsion Free Sheaves}
\label{lena}
\subsection{Step V}
Let  $g \geq2$, $r>0$ be integers. Let
$$e > max(e_1(g,r), e_2(g,r)),$$
$$t> max(t_0(g,r,e), t_1(g,r,e), t_2(g,r,e)),$$
be determined by
Propositions (\ref{wk1}, \ref{bob}, \ref{kiwi}, \ref{cake}).
We now conclude  the only
possible semistable points in the $SL_n$-linearized
G.I.T. problem determined
by
$$i_t: Q_g(C,n=f_{e,r}(0), f_{e,r}) \rarr
\grass(f_{e,r}(t), (\tp Sym^t(H^0(C,w_C^{10})))^*) $$
are elements $\xi\in Q_g(C,n,f_{e,r})$ that
correspond to quotients
$$\tp \oh_C \rarr E \rarr 0$$
where $$\psi:\tp H^0(C,\oh_C) \rarr H^0(C,E)$$
is an isomorphism and $E$ is a slope-semistable, torsion
free sheaf.
In order for $\xi$ to be semistable,
$\psi$ must be injective by Proposition (\ref{wk1}).
Surjectivity is equivalent to $h^1(C,E)=0$. By Proposition
(\ref{cake}), $\psi$ must be surjective. Since $\psi$ is
an isomorphism, $E$ must be torsion free by Proposition (\ref{bob}).
Finally, by Proposition (\ref{kiwi}), $E$ must be slope-semistable.
We now establish the converse.

\begin{pr}
\label{snake}
Let $g\geq 2$, $r>0$ be integers. There exist bounds
$e_3(g,r)>r(g-1)$ and $t_3(g,r,e)>\hat{t}(g,r,e)$
such that for each pair $e>e_3(g,r)$, $t>t_3(g,r,e)$ and any
$[C]\in H_g$, the following holds:

\noindent
If $\ \xi \in Q_g(C,n=f_{e,r}(0),f_{e,r})$ corresponds to
a quotient
$$\tp \oh_C \rarr E \rarr 0$$
where
$$ \psi: \tp H^0(C,\oh_C) \rarr H^0(C,E)$$
is an isomorphism and $E$ is a slope-stable (slope-semistable),
torsion free sheaf,
then $\xi$ is a G.I.T. stable (semistable) point for the
$SL_n$-linearization determined by $i_t$.
\end{pr}

\subsection{Lemmas}
For the proof of G.I.T. (semi)stability, the fundamental step
is the inequality of Lemma (\ref{king}) for every subsheaf $F$ of
$E$ generated by global sections.
Note this is the reverse of the inequality required by
Lemma (\ref{archy}). Lemma (\ref{king}) follows by vanishing, Riemann-Roch,
and the slope-(semi)stability of $E$ when $F$ is nonspecial. In case
$h^1(C,F)\neq 0$, an analysis in Lemma (\ref{uinta})
utilizing $Hom(F,\om)\neq 0$ yields the
required additional information. The Numerical Criterion of
section (\ref{nnmc}) and Lemma (\ref{king})
reduce the stability question to a purely
combinatorial result established in Lemma (\ref{queen}).

\begin{lm}
\label{uinta}
Let $g\geq 2$, $r>0$ be integers. Let $q$ be an integer.
There exists an integer
$e_3(g,r,q)$ such that for each
$e >e_3(g,r,q)$ and $[C]\in H_g$, the following holds:

\noindent
If $E$ is  a slope-semistable, torsion free sheaf on
$C$ with Hilbert polynomial $f_{e,r}$ (with respect
to $\om_C^{10}$) and
$$0 \rarr F \rarr E$$ is  a nonzero subsheaf with multirank
$(s_i)$ satisfying
$h^1(C,F) \neq 0$, then:
\begin{equation}
\label{hot}
{\chi(F) +q \over \sum s_i\om_i}
< {\chi(E) \over r(2g-2)}-1.
\end{equation}
\end{lm}

\begin{pf}
Since $h^1(C,F) \neq 0$, there exists a nontrivial morphism
$$\sigma: F \rarr \om_C.$$
Consider the subsheaf
$0 \rarr \sigma(F) \rarr \om_C $
where $\sigma(F) \neq 0$.  By the proof of Lemma (\ref{boy}),
$\chi(\sigma(F)) \leq g$. Consider the exact sequence
$$0\rarr K \rarr F \rarr \sigma(F) \rarr 0.$$
If $K=0$, then
$$ {\chi(F) +q \over \sum s_i\om_i}
< g+q. $$
Therefore, if
$$e> r(g-1) + r(2g-2)(g+q+1),$$ the case $K=0$ is
settled. Also,  the case $\chi(F) \leq g$ is settled.
Now suppose $K\neq 0$ and $\chi(F)-g >0$.
We have $\chi(F)-g  \leq \chi(K)$.
 Let $(s'_i)$ be the nontrivial
multirank of $K$.
Since $\sigma(F)$ is of nontrivial multirank we have:
$$ 0 < \sum s'_i \om_i < \sum s_i \om_i \leq r(2g-2).$$
We obtain:
\begin{eqnarray*}
\paren{\chi(F) -g \over \sum s_i \om_i} \cdot
 \paren{r(2g-2) \over r(2g-2)-1}
& \leq & \paren{\chi(F) -g \over \sum s_i \om_i } \cdot
 \paren{\sum s_i\om_i \over \sum s'_i\om_i} \\
& \leq & {\chi(K) \over \sum s'_i \om_i} \ \  .
\end{eqnarray*}
Using the slope-semistability of $E$ with respect to $K$,
we conclude:
$${\chi(F) -g \over \sum s_i \om_i}\leq
\paren{\chi(E) \over  r(2g-2)}
 \cdot
 \paren{r(2g-2)-1 \over r(2g-2)} \ .$$
It is now clear, for large $e$ depending
only on $g$ and $r$ and $q$, the inequality (\ref{hot})
is satisfied.
\end{pf}

\begin{lm}
\label{king}
Let $g\geq2$, $r> 0$. Let $b=b(g,R(g,r))$.
Let $e>e_3(g,r,b)>r(g-1)$.
There exists an integer
 $t_3(g,r,e)> \hat{t}(g,r,e) $
such that for each
$t > t_3(g,r,e)$ and $[C]\in H_g$, the following holds:

\noindent
If $\xi \in Q_g(C,n=f_{e,r}(0), f_{e,r})$
corresponds to a quotient
$$ \tp \oh_C \rarr E \rarr 0$$
where $E$ is a torsion free, slope-semistable
sheaf on $C$ and
$0 \rarr F \rarr E$ is a nonzero, proper subsheaf generated by
global sections, then
$$ {h^0(C,F)\over n}  \leq
{h^0(C,F\otimes \om_C^{10t})\over f_{e,r}(t)}\ .$$
If E is slope-stable,
$$ {h^0(C,F)\over n}  <
{h^0(C,F\otimes \om_C^{10t})\over f_{e,r}(t)}\ .$$
\end{lm}

\begin{pf}
Suppose $t> p(g, R(g,r), \chi=f_{e,r}(0))$.
Let $(s_i)$ be the nontrivial multirank of
$F$. Since $F$ is generated by global sections,
the vanishing
guaranteed by Lemma (\ref{cotton}) yields
$$h^0(C,F\otimes \om_C^{10t}) = \chi(F) + (\sum s_i\om_i)10t.$$
The Hilbert polynomial can be expressed:
$$f_{e,r}(t) = \chi(E) + r(2g-2)10t.$$
First consider the case where $h^1(C,F)=0.$
Then $h^0(C,F) =\chi(F)$.
We compute
$$ \chi(F) \cdot f_{e,r}(t)-
 \chi(E) \cdot h^0(C,F\otimes \om_C^{10t}) =$$
$$\chi(F)\cdot r(2g-2)10t - \chi(E) \cdot (\sum s_i\om_i)10t \ < (\leq) \  0.$$
where $E$ is slope-stable, (slope-semistable). Hence
$$ {h^0(C,F)\over \chi(E)} =
 {\chi(F)\over \chi(E)}  < (\leq) \
{h^0(C,F\otimes \om_C^{10t})\over f_{e,r}(t)}\ .$$
Since $n=\chi(E)$,
the nonspecial case is thus settled.
Now suppose $h^1(C,F)\neq0$.
Lemma (\ref{uinta}) now applies to $F$.
We compute
$$ (\chi(F)+b) \cdot f_{e,r}(t)-
 \chi(E) \cdot h^0(C,F\otimes \om_C^{10t}) =$$
$$(\chi(F)+b)\cdot r(2g-2)10t - \chi(E) \cdot (\sum s_i\om_i)10t
+ b \cdot \chi(E) <$$
$$ -(\sum s_i \om_i) \cdot r(2g-2) \cdot 10 t + b\cdot \chi(E).$$
For $t > b \cdot \chi(E) = b\cdot (e+r(1-g))$,
$$ {\chi(F)+b\over  \chi(E)}  <
{h^0(C,F\otimes \om_C^{10t})\over f_{e,r}(t)}\ .$$
Since $h^0(C,F)<\chi(F)+b$,
$${h^0(C,F)\over \chi(E)} < {\chi(F)+b \over  \chi(E)}$$
The proof is complete.
\end{pf}

We require a simple combinatorial Lemma.
\begin{lm}
\label{queen}
Let $n \geq 2$ be an integer.
Let $$W_1 \leq W_2 \leq \ldots \leq W_n , \ \ W_1 <W_n$$
be integers. Let $\{\al_i\}$, $ \{ \bay_i \} $ be rational numbers
such that
\begin{enumerate}
\item[(i.)] $\sum_{1}^{n}\bay_i= \sum_{1}^{n}\al_i$.
\item[(ii.)]
 $\forall \ 1 \leq m \leq n-1$, $\ \ \sum_{1}^{m}\bay_i \  <  (\leq)\
\sum_{1}^{m}\al_i$.
\end{enumerate}
Then:
$$\sum_{1}^{n}\bay_i \cdot W_i\  > (\geq) \  \sum_{1}^{n}\al_i \cdot W_i.$$
\end{lm}
\begin{pf}
Use discrete Abel summation:
$$\sum_{i=1}^{n}\bay_i \cdot W_i = (\sum_{i=1}^{n}\bay_i)\cdot W_n  -
\sum_{m=1}^{n-1} \Bigg( (\sum_{1}^{m}\bay_i)\cdot (W_{m+1}-W_m) \Bigg)
 \ \ > (\geq)$$
$$(\sum_{i=1}^{n}\al_i)\cdot W_n  -
\sum_{m=1}^{n-1} \Bigg( (\sum_{1}^{m}\al_i)\cdot (W_{m+1}-W_m) \Bigg) =
\sum_{i=1}^{n}\al_i \cdot W_i.$$
The middle inequality follows from (i) and (ii) above.
\end{pf}

\subsection{Proof of Proposition (\ref{snake})}

\begin{pf}
Let $e_3(g,r)=e_3(g,r,b(g,R(g,r)))$
be determined  by Lemma (\ref{uinta}). For $e>e_3(g,r)$, let
$t_3(g,r,e)> p(g,R(g,r), \chi=f_{e,r}(0))$
 be given by Lemmas (\ref{king}) and (\ref{cotton}).
We will apply the Numerical Criterion to the linearized
$SL_n$-action on
$$\grass(f_{e,r}(t), (\bold{C}^n \otimes Sym^t(H^0(C,\om_C^{10})))^*).$$
The element $\xi$ corresponds to the quotient:
$$\psi^t:\bold{C}^n \otimes Sym^t(H^0(C,\om_C^{10})) \rarr
H^0(C,E\otimes \om_C^{10t}) \rarr 0.$$
Let $\barr{v}=(v_1, \ldots, v_n)$ be a basis of $\bold{C}^n$.
Let $(w(v_1),\ldots, w(v_n))$ be  weights satisfying
$$ w(v_1) \leq w(v_2) \leq \ldots \leq w(v_n), \ \ w(v_1) < w(v_n).$$
To apply the Numerical Criterion for (semi)stability, an $f_{e,r}(t)$-tuple
of $\barr{v}$-pure
elements of $\tp Sym^t(H^0(C, \om_C^{10}))$ projecting to a basis
of $H^0(C, E\otimes \om_C^{10})$ and satisfying the weight inequality (2)
of section (\ref{nnmc}) must be shown to exist.

For $1\leq i \leq n$,
let $F_i$ denote the subsheaf of $E$ generated by
$\psi(\sumo_{j=1}^{i} v_j\otimes H^0(C,\oh_C)).$
By the surjectivity guaranteed by (ii) of Lemma (\ref{cotton}),
\begin{equation}
\label{brain}
\psi^t : \sumo_{j=1}^{i} v_j\otimes Sym^t(C,H^0(\om_C^{10})) \rarr
H^0(C,F_i\otimes \om_C^{10t}) \rarr 0.
\end{equation}
Define for  $1\leq i \leq n$, $\ A_i= h^0(C,F_i\otimes \om_C^{10t})$.

The required $f_{e,r}(t)$-tuple
$(a_1, a_2, \ldots, a_{f_{e,r}(t)})$ is constructed as
follows. Select  elements $$(a_1, \ldots, a_{A_1})\in \
v_1 \otimes Sym^t(H^0(C,\om_C^{10}))$$ such that
$(\psi^t(a_1), \ldots, \psi^t(a_{A_1}))$
determines a basis of $H^0(C,F_1\otimes \om_C^{10t})$.
Select $$(a_{A_1+1}, \ldots , a_{A_2}) \in \
v_2 \otimes Sym^t(H^0(C,\om_C^{10}))$$ such that
$(\psi^t(a_1), \ldots, \psi^t(a_{A_2}))$ determines  a basis
of $H^0(C,F_2\otimes \om_C^{10t})$. Continue selecting
$$(a_{A_i+1}, \ldots , a_{A_{i+1}}) \in \
v_{i+1} \otimes Sym^t(H^0(C,\om_C^{10}))$$ such that
$(\psi^t(a_1), \ldots, \psi^t(a_{A_{i+1}}))$ determines a basis
of $H^0(C,F_{i+1}\otimes \om_C^{10t})$. Note if $A_{i}=A_{i+1}$,
then $(\psi^t(a_1), \ldots, \psi^t(a_{A_i}))$ already
determines a basis of $H^0(C,F_{i+1} \otimes \om_C^{10t})$ and
no elements of $v_{i+1}\otimes Sym^t(H^0(C,\om_C^{10}))$
are chosen.
This selection
is possible by the surjectivity of (\ref{brain}).

Let $\al_1= A_1/ f_{e,r}(t)$ and $\al_i= (A_i-A_{i-1})/f_{e,r}(t)$ for
$2\leq i \leq n$. We have
$$\sum_{i=1}^{n}\al_iw(v_i) = \sum_{j=1}^{f_{e,r}(t)}{w(a_j)\over
f_{e,r}(t)}.$$
Let $\bay_i= (1/n)$.
Note $$\sum_{1}^{n}\bay_i=\sum_{1}^{n}\al_i  = 1.$$
Let $1 \leq m \leq n-1$.
Since $\psi$ is an isomorphism, $m\leq h^0(C,F_m)$.
Suppose $F_m \neq E$.
Then Lemma (\ref{king}) yields
\begin{equation}
\label{dayo}
{m\over n} \ < (\leq) \  {A_m\over f_{e,r}(t)}\ .
\end{equation}
If $F_m$= $E$, then $A_m=f_{e,r}(t)$ and
 the inequality (\ref{dayo}) holds trivially ($m\leq n-1$).
So for all $1\leq m \leq n-1$,
we have $$ \sum_{1}^{m}\bay_i\  < (\leq)\  \sum_{1}^{m}\al_i.$$
Lemma (\ref{queen}) yields
$$\sum_{i=1}^{n}{w(v_i)\over n} =
\sum_{1}^{n} \bay_i w(v_i) \ > (\geq)\  \sum_{1}^{n}  \al_i w(v_i)
 = \sum_{j=1}^{f_{e,r}(t)}{w(a_j)\over f_{e,r}(t)} .$$
By the Numerical Criterion, $\xi$ is G.I.T. stable (semistable).
\end{pf}

\subsection{Step VI}
Only one step remains in the proof of Theorem (\ref{fred}).
It must be checked that strict slope-semistability of $E$ implies
strict G.I.T. semistability.
\begin{lm}
\label{moon}
Let $g\geq 2$, $r>0$ be integers.
There exists an integer $e_4(g,r)$ such that for each
$e>e_4(g,r)$ and $[C]\in H_g$, the following holds:

\noindent
If $E$ is any  slope-semistable, torsion free sheaf on a
$C$
with  Hilbert polynomial $f_{e,r}$
(with respect to $\om_C^{10t}$) and
$0 \rarr F \rarr E$ is
a nonzero subsheaf with multirank $(s_i)$ satisfying
\begin{equation}
\label{tart}
{\chi(F)\over \sum s_i \om_i}=
{\chi(E) \over r(2g-2)},
\end{equation}
then
\begin{enumerate}
\item[(i.)] $h^1(C,F)=0.$
\item[(ii.)] $F$ is generated by global sections.
\end{enumerate}
\end{lm}

\begin{pf}
Suppose $F$ is a nonzero subsheaf of $E$ satisfying (\ref{tart}).
If $e>e_3(g,r,0)$, by Lemma (\ref{uinta}), $h^1(C,F)=0$.
Now let $x\in C$ be a point. We have an exact sequence
$$0 \rarr m_xF \rarr F \rarr {F\over m_xF} \rarr 0.$$
Since $F$ is torsion free and $C$ is nodal, it is not hard to show that
$$dim \paren{F\over m_xF} < 2\cdot R(g,r).$$
Since $(F/m_xF)$ is torsion, $m_xF$ has the same  multirank
as $F$.  Also $$\chi(m_xF) > \chi(F) - 2R.$$
By (\ref{tart}),
$${\chi(m_xF)+2R \over \sum s_i \om_i} >
{\chi(E)\over  r(2g-2)}. $$
If $e> e_3(g,r, 2R)$, $h^1(C,m_xF)=0$ by Lemma (\ref{uinta}).
In this case $F$ is generated by global sections.
We can therefore choose
$e_4(g,r) = e_3(g,r,2R)$.
\end{pf}

\begin{pr}
\label{snook}
Let $g\geq 2$, $r>0$ be integers. There exist bounds
 $e_4(g,r)$ and $t_4(g,r,e)$ such that
for each pair $e>e_4(g,r)$, $t>t_4(g,r,e)$ and any
$[C]\in H_g$, the following holds:

\noindent
If $\xi \in Q_g(C,n=f_{e,r}(0),f_{e,r})$ corresponds to
a quotient
$$\tp \oh_C \rarr E \rarr 0$$
where
$$ \psi: \tp H^0(C,\oh_C) \rarr H^0(C,E)$$
is an isomorphism and $E$ is a strictly slope-semistable,
torsion free sheaf,
then $\xi$ is a G.I.T. strictly semistable point for the
$SL_n$-linearization determined by $i_t$.
\end{pr}

\begin{pf}
Since $\xi$ is G.I.T. semistable for $e>e_3(g,r)$ and
$t>t_3(g,r,e)$, it suffices to find a semistabilizing
$1$-parameter subgroup. If $e>e_4(g,r)$, then, for  any nonzero, proper
semistabilizing subsheaf
$0\rarr F \rarr E$,
we have $h^0(C,F)= \chi(F)$ and $F$ is generated
by global sections. It is now easy to see that
the flag $0\subset H^0(C,F) \subset H^0(C,E)$ with
weights $\{0,1 \}$ determines semistabilizing data for large
$t>t_4(g,r,e)$.
\end{pf}

We have now shown for  the bounds:
$$ e(g,r) = max \{ e_i(g,r) \ | \ 1 \leq i \leq 4 \},$$
$$t(g,r,e) = max \{t_i(g,r,e) \ | \ 0 \leq i \leq 4 \}$$
the claim of Theorem (\ref{fred}) holds. This completes the proof of
Theorem (\ref{fred}).

By Lemma (\ref{moon}),   each slope-semistable, torsion free
sheaf $E$ on $C$ with Hilbert polynomial $f_{e,r}$
 appears in $Q_g(C,f_{e,r}(0),f_{e,r})^{SS}_{t}$ for $e>e(g,r)$,
$t>t(g,r,e)$. It is now
clear the $SL_n$-orbits of $Q_g(C,f_{e,r}(0), f_{e,r})^{SS}_{t}$
correspond exactly to the slope-semistable, torsion free
sheaves on $C$ with Hilbert polynomial $f_{e,r}$.

\subsection{Seshadri's Construction}
\label{sesh}
In [Se],
C. Seshadri has studied the $SL_n$-action on $Q_g(C,n=f_{e,r}(0),f_{e,r})$
via a covariant construction. For $e>>0$, he finds a G.I.T.
(semi)stable locus that coincides exactly with the
G.I.T. (semi)stable locus of Theorem (\ref{fred}). These results
appear in Theorem 18 of chapter 1 of [Se] for nonsingular curves and
Theorem 16 of chapter 6 for singular curves.
 The
collapsing of semistable orbits is determined by the
Zariski topology. Seshadri shows that
\begin{enumerate}
\item[(i.)]  If
$E_{t}$ is a  flat family of
of slope-semistable, torsion free sheaves on $C$ such that the
Jordan-Holder factors of $E_{t\neq 0}$ are $\{ F_j \}$, then
the Jordan-Holder factors of $E_0$ are also $\{ F_j \}$.
\item[(ii.)] If $E$ is a
slope-semistable, torsion free sheaf on $C$ with
Jordan-Holder factors $\{ F_j \}$, then there exists a flat
family of slope-semistable, torsion free sheaves $E_t$
such that:
 $$E_{t\neq 0} \cong E, \ \ \  E_0  \cong \bigoplus_{j} F_j\ \ .$$
\end{enumerate}
Statement (ii) is proven by constructing flat families
over extension groups.
It follows from these two results that
the points of our  quotient
$$Q_g(C,n=f_{e,r}(0),f_{e,r})_{t}^{SS} /SL_n$$
for $e>e(g,r),\  t>t(g,r,e)$
naturally parametrize slope-semistable, torsion free sheaves
with Hilbert polynomial $f_{e,r}$ up to equivalence
given by Jordan-Holder factors.

\section{Two Results in Geometric Invariant Theory}
\label{abgit}
\subsection{Statements}
Let $V$, $Z$, and $W$ be finite dimensional $\bold{C}$-vector
spaces.
Consider two rational representations of $SL(V)$:
$$\zeta: SL(V) \rarr SL(Z) $$
$$\omega: SL(V) \rarr SL(W). $$
These representations define natural $SL(V)$-linearized actions
on $\proj(Z)$ and $\proj(W)$. There is an induced $SL(V)$-action
on the product $\proj(Z) \times \proj(W)$. Since
$$Pic(\proj(Z) \times \proj(W)) = \bold{Z} \oplus \bold{Z},$$
there is a  $1$-parameter choice of linearization.  For
$a,b \in \bold{N}^+$, let $[a,b]$ denote the linearization
given by the line bundle $\oh_{\proj(Z)}(a) \otimes
\oh_{\proj(W)}(b)$.  Subscripts will be used to indicate
linearization. Let
$$\rho_Z : \proj(Z) \times \proj(W) \rarr \proj(Z)$$
be the projection on the first factor.

\begin{pr}
\label{jack}
There exists an integer $k_S(\zeta,\om)$ such that for all $k>k_S$:
$$\rho_Z^{-1}(\proj(Z)^S) \subset (\proj(Z) \times \proj(W))_{[k,1]}^S.$$
\end{pr}

\begin{pr}
\label{jill}
There exists an integer $k_{SS}(\zeta,\om)$ such that for all $k>k_{SS}$:
$$(\proj(Z) \times \proj(W))_{[k,1]}^{SS} \subset
\rho_Z^{-1}(\proj(Z)^{SS}).$$
\end{pr}

D. Edidin has informed the author that Proposition (\ref{jack})
is essentially equivalent to Theorem 2.18 of [M-F].

\subsection{q-Stability}
\label{nmc}
Let
$\lambda: \bold{C}^* \rarr SL(V)$ be a $1$-parameter subgroup.
Let $dim(V)=a$.
It is well known there exists  a basis $\barr{v}=(v_1, \ldots, v_a)$
of $V$ such that $\lambda$ takes the form
$$\lambda(t)(v_i)=t^{e_i} \cdot v_i, \ \ \ \ t\in \bold{C}^*.$$
Denote the tuple $(e_1, \ldots, e_a)$ by $\barr{e}$.
The exponents satisfy the determinant $1$ condition,
$\sum_{i=1}^{a}e_i=0$. Let $|\barr{e}|= max\{ |e_i| \}$.
For the representation $\zeta:SL(V) \rarr SL(Z)$, there
exists a basis $\barr{z}=(z_1, \ldots, z_b)$ such that
$\zeta \circ \lambda$ takes the form
$$\zeta \circ \lambda(t)(z_j)=t^{f_j} \cdot z_j,\ \ \ \ t\in \bold{C}^*.$$
Again, $\sum_{j=1}^{b} f_j=0.$ The pairs $\{\barr{v}, \barr{e} \}$ and
$\{\barr{z}, \barr{f} \}$ are said to be {\em diagonalizing data}
for $\lambda$ and $\zeta \circ \lambda$ respectively.

 Let
$[z]\in \proj(Z)$ correspond to the one dimensional subspace
of $Z$ spanned by $z\neq 0$. By the Mumford-Hilbert Numerical
Criterion, $[z]$ is a stable (semistable) point
for the $\zeta$-induced
linearization
on $\proj (Z)$ if and only if for every $1$-parameter subgroup
$\lambda: \bold{C}^* \rarr SL(V)$, the following condition
holds: let $\{\barr{z}, \barr{f} \}$ be diagonalizing data for
$\zeta \circ \lambda$ and let $z=\sum_{j=1}^{b} \xi_j \cdot z_j$,
then there exists an index  $j$ for which $\xi_j \neq 0$ and $f_j <0$
($f_j \leq 0$).

Let $q>0$ be a real number.
The point $[z]$ is defined to be {\em $q$-stable} for the
$\zeta$-induced linearization  if and only if for  every $1$-parameter
 subgroup
$\lambda: \bold{C}^* \rarr SL(V)$ the following condition
holds: let $\{\barr{v}, \barr{e} \}$ and $\{\barr{z}, \barr{f} \}$
 be diagonalizing data for $\lambda$ and
$\zeta \circ \lambda$ and let $z=\sum_{j=1}^{b} \xi_j \cdot z_j$,
then there exists an index  $j$ for which $\xi_j \neq 0$ and
$f_j < -q \cdot  |\barr{e}|.$
For $q>0$, let $\proj(Z)^{qS}$ denote the $q$-stable locus
for the $\zeta$-induced linearization.

Proposition (\ref{jack}) will be established in two
steps:
\begin{lm}
\label{step1}
There exists $q(\zeta)>0$ such that $\proj(Z)^{qS}=\proj(Z)^S$.
\end{lm}

\begin{lm}
\label{step2}
For any $q>0$,
there exists an integer $k_{qS}(q,\om)$ such that for all $k>k_{qS}$:
$$\rho_Z^{-1}(\proj(Z)^{qS}) \subset (\proj(Z) \times \proj(W))_{[k,1]}^S.$$
\end{lm}
\noindent
Lemmas (\ref{step1}) and (\ref{step2}) certainly imply
Proposition (\ref{jack}).

\subsection{Proofs of Lemmas (\ref{step1}) and (\ref{step2})}
Let $U$ be a finite dimensional $\bold{Q}$-vector space.
Let $L=\{l_i \}$ be a finite set of elements of $U^*$.
The set $L$ is said to be  a {\em stable configuration} if
$$\forall\  0\neq u\in U,  \ \ \exists i \ \ l_i(u)<0.$$
If $\barr{u}=(u_1,u_2,\ldots,u_d)$ is a basis of $U$,  define
a norm $|\ |_{\barr{u}}: U\rarr \bold{Q}^{\geq 0}$ by
$$|u|_{\barr{u}}= max \{ |\gamma_i|  \}\ \ \ \mbox{where }\ u=
\sum_{1}^{k}\gamma_iu_i .$$

\begin{lm}
\label{sam}
Suppose $L =\{l_i \}$ is a stable configuration in $U$.
Let $\barr{u}$ be a basis of $U$.
Then there exists $q>0$ depending upon $L$ and $\barr{u}$ such that
\begin{equation}
\label{drop}
\forall \  0\neq u\in U, \ \ \exists i \ \ l_i(u)<-q\cdot |u|_{\barr{u}}.
\end{equation}
\end{lm}

\begin{pf}
Let $U \subset U_{\bold {R}} \cong U \otimes_{\bold{Q}} \bold{R} $.
Suppose there exists an element  $0\neq u\in U_{\bold{R}}$ and a decomposition
$L=L' \cup L''$  satisfying:
\begin{enumerate}
\item[(i.)] $ \forall l\in L',  \ \  l(u)=0.$
\item[(ii.)] $\forall l\in L'',  \ \l(u)>0.$
\end{enumerate}
Since the locus $\{z\in U_{\bold{R}}\  | \ \forall l\in L',\   l(z)=0 \}$ is a
rational subspace and the locus
 $\{ z\in U_{\bold{R}}\  | \ \forall l\in L'', \  l(z)>0 \}$ is
open, there must exist an element  $0\neq \hat{u}\in U$ satisfying
(i) and (ii).
Since $L$ is a stable configuration in  $U$, such
$\hat{u}$ do not exist.
It follows
$$\forall \  0\neq u\in U_{\bold{R}}, \ \ \exists i \ \ l_i(u)< 0.$$
Let $S$ be the unit box in $U_{\bold{R}}$:
$S= \{ u\in U_{\bold{R}}\ | \ |u|_{\barr{u}}=1\}.$
Define a function $g: S \rarr \bold{R}^-$ by
$$g(s)= min\  \{l(s) \ |\  l\in L \}.$$
The function $g$ is continuous and strictly negative. Since
$S$ is compact, $g$ achieves a maximum value $-m$ on $S$ for some
$m>0$. The bound $q=m/2$ clearly satisfies (\ref{drop}).
\end{pf}

\begin{pf}[Of Lemma (\ref{step1})]
The proof consists of two simple pieces.
First, a basis $\barr{v}$ of $V$ is fixed. By applying Lemma (\ref{sam}),
it is shown there exists a $q>0$ such that the
stability of $[z]$ implies the $q$-stability inequality for all
$1$-parameter subgroups
subgroups of $SL(V)$ diagonal with respect $\barr{v}$.
Second, it is checked that this $q$ suffices for any selection of
basis.

Let $\barr{v}=(v_1, \ldots, v_a)$ be a basis of $V$.
Let $$U= \{(e_1,\ldots,e_a) \ |\ \  e_i\in \bold{Q}, \sum_{1}^{a}e_i=0 \}.$$
There exist linear functions $\{ l_1, \ldots , l_b \}$ on $U$ and
a basis $\barr{z}=(z_1, \ldots, z_b)$ of $Z$ satisfying the following:
if $\lambda : \bold{C}^* \rarr SL(V)$ is any $1$-parameter subgroup
with diagonalizing data $(\barr{v}, \barr{e})$, then the
diagonalizing data of $\zeta \circ \lambda$ is
$(\barr{z}, (l_1(\barr{e}), \ldots, l_b(\barr{e})))$.
Let $\{L_1, \ldots, L_{B} \}$ be the set of distinct
stable configurations in $\{l_1, \ldots ,l_b \}$.  That is,
for all $1\leq J  \leq B $,
 $\ L_{J} \subset \{l_1, \ldots ,l_b \}$ and $L_{J}$ is
a stable configuration in $U$. Let $\barr{u}=(u_1, \ldots, u_{a-1})$
be a basis of $U$ of the following form:
 $$u_1=(-1, 1,0,\ldots,0),\  \ldots\ ,\  u_{a-1}=(-1,0,\ldots,0,1).$$
Note $|\barr{e}| \leq a \cdot |\barr{e}|_{\barr{u}}$ for $\barr{e}\in U$.
By Lemma (\ref{sam}), there exists
$q_{J}>0$ such that (\ref{drop}) holds
for each stable configuration $L_{J}$.
Let $$q={1\over a} \cdot min\{q_{J} \}.$$
Suppose $[z]\in \proj(Z)^S$. Let $z=\sum_{1}^{b} \xi_j z_j$.
By the Numerical Criterion, the stability of $[z]$ implies the
set $\{ l_j | \xi_j \neq 0 \}$ is a stable
configuration in $U$ equal to some $L_{J}$. For any $1$-parameter
subgroup with diagonalizing data
$({\barr{v}}, {\barr{e}})$,
the  diagonalizing data of $\zeta \circ \lambda$ is
$(\barr{z}, (l_1(\barr{e}), \ldots, l_b(\barr{e})))$.
By (\ref{drop}), we see
there exists an $l_i\in L_{J}$ such that
$$ l_i(\barr{e})< -q_{J} \cdot |\barr{e}|_{\barr{u}}
 \leq -q \cdot |\barr{e}|.$$
Suppose $\barr{v}'$ is another basis of $V$. Then, up to scalars, there
exists an element $\gamma \in SL(V)$ satisfying $\gamma (\barr{v})=
\barr{v}'$. It is now clear that
$$(\zeta(\gamma)(\barr{z}), (l_1(\barr{e}), \ldots, l_b(\barr{e})))$$
 is diagonalizing  data for $\zeta \circ \lambda$ where
$\lambda$ has diagonalizing data $(\barr{v}',\barr{e})$.
Since the set $\{ l_1, \ldots, l_b \}$ is independent of $\barr{v}$,
the above analysis is valid for any $1$-parameter subgroup.
We have shown that $[z]\in \proj(Z)^{qS}.$
\end{pf}

\begin{lm}
\label{wally}
Let $\omega: SL(V) \rarr SL(W)$ be a rational representation. There
exists an $M_\omega > 0$ with the following property.
Let $\lambda: \bold{C}^* \rarr SL(V)$ be any  $1$-parameter
subgroup. Let $(\barr{v}, \barr{e})$ and $(\barr{w},\barr{h})$
be diagonalizing data for $\lambda$ and $\omega \circ \lambda$.
Then $|h| < M_\omega \cdot |e|$.
\end{lm}

\begin{pf}
Let $\barr{v}$ be a basis of $V$. Let $U$ be as in the
proof of Lemma (\ref{step1}). There exist
linear functions $\{ l_1, \ldots, l_c \}$ on $U$  and a basis
$\barr{w}=(w_1, \ldots, w_c)$ of
$W$ satisfying the following: if $\lambda: \bold{C}^* \rarr
SL(V)$ is any 1-parameter subgroup with diagonalizing data
$(\barr{v}, \barr{e})$, then the diagonalizing data
of $\omega \circ \lambda$ is $(\barr{w}, (l_1(\barr{e}),
\ldots, l_c(\barr{e})))$. Choose $M_\omega$ so
$$\forall j, \ \ |l_j(\barr{e})| < M_\omega \cdot |\barr{e}|.$$
As in the proof of Lemma (\ref{step1}), the set of linear
functions does not depend on $\barr{v}$.  The proof is complete.
\end{pf}

\begin{pf} [Of Lemma (\ref{step2})] It is clear that if an
element $[z]\in\proj(Z)$ is $q$-stable for the $\zeta$-induced
linearization, then $[z^k]\in \proj(Sym^k(Z))$ is
$kq$-stable for the $Sym^k(\zeta)$-induced
linearization. Let $M_\omega$ be determined by  Lemma (\ref{wally})
for the representation $\omega$. Let $k_{qS}=M_\omega/q$.
We check for $k>k_{qS}$,
$$\rho_Z^{-1}(\proj(Z)^{qS}) \subset (\proj(Z) \times \proj(W))_{[k,1]}^S.$$
The linearization $[k,1]$ corresponds to the embedding:
$$\proj(Z) \times \proj(W) \rarr \proj(Sym^k(Z) \otimes W)$$
$$[z]\times [w] \rarr [z^k\otimes w].$$
Let $[z]\in \proj(Z)^{qS}$ and $[w] \in \proj(W)$. Let
$\lambda: \bold{C}^* \rarr SL(V)$ be any $1$-parameter subgroup.
Let $\barr{e}$ be the diagonalized exponents of $\lambda$.
Let $(\barr{z}^*, \barr{f}^*)$ and $(\barr{w}, \barr{h})$
be the diagonalizing data of $Sym^k(\zeta) \circ \lambda$ and
$\omega \circ \lambda$. Since $[z^k]$ is $kq$-stable for the
$Sym^k(\zeta)$-induced linearization, there exists an index
$\mu$ satisfying:
\begin{enumerate}
\item [(i.)] The basis element $z_\mu^*$ has a nonzero coefficient
in the expansion of $[z^k]$.
\item [(ii.)] $f_\mu^* < - kq \cdot |\barr{e}|
< -M_\omega \cdot |\barr{e}|$.
\end{enumerate}
Let $w_\nu$ be any basis element that has a nonzero
coefficient in the expansion of $w$. Note
$z_\mu^* \otimes w_\nu$ is an element of the
diagonalizing basis $\barr{z}^* \otimes \barr{w}$ of
$$(Sym^k(\zeta)\otimes \omega) \circ \lambda$$ having
nonzero coefficient in the expansion of $z^k\otimes w$.
The exponent corresponding to $z_\mu^* \otimes w_\nu$ is
simply $f_\mu^* + h_\nu$. Since
$$|h_\nu| \leq |\barr{h}| < M_\omega \cdot |e|,$$
condition (ii) above implies the exponent is
strictly negative.  By the Numerical Criterion,
$[z^k\times w]$ is stable. The Lemma is proven.
\end{pf}

\subsection{Proof of Proposition (\ref{jill})}
Let $\zeta: SL(V) \rarr SL(Z)$
be a rational representation as above. An element $[z]\in \proj(Z)$ is
{\em $(e_1, \ldots, e_a)$-unstable} for the $\zeta$-induced
linearization if there exists a destabilizing $1$-parameter subgroup
$\lambda: \bold{C}^* \rarr SL(V)$ with diagonalizing data
$(\barr{v}, \barr{e})$: if $(\barr{z}, \barr{f})$ is
diagonalizing data for $\zeta \circ \lambda$ and
$z= \sum_{1}^{b} \xi_j \cdot z_j$, then
$\xi_j\neq 0$ implies $f_j >0.$
Let $\proj(Z)^{\barr{e}UN} \subset \proj(Z)$ denote
$\barr{e}$-unstable locus.

{}From the Numerical Criterion, every unstable point is
$\barr{e}$-unstable for some $a$-tuple $\barr{e}=(e_1,\ldots, e_a)$. We need
a simple finiteness result:

\begin{lm}
\label{tom}
Consider the $\zeta$-linearized G.I.T. problem on $\proj(Z)$.
There exists a finite set of
$a$-tuples, $\cal{P}$,
such that
$$\bigcup_{\barr{e}\in \cal{P}} \proj(Z)^{\barr{e}UN} = \proj(Z)^{UN}.$$
\end{lm}

\begin{pf}
We first show that
$\proj(Z)^{\barr{e}UN}$ is a constructible
subset of $\proj(Z)$. Fix a $1$-parameter subgroup
$\lambda: \bold{C}^* \rarr SL(V)$ with diagonalizing
data $(\barr{v}, \barr{e})$.  Let $(\barr{z}, \barr{f})$
be diagonalizing data for $\zeta \circ \lambda$. Let $H$
be the projective subspace of $\proj(Z)$ spanned by the set
$\{z_j| f_j>0 \}$. Certainly $H\subset  \proj(Z)^{\barr{e}UN}$.
Since every $1$-parameter subgroup of $SL(V)$ with
diagonalized exponents $\barr{e}$ is conjugate to
$\lambda$, we see the map:
$$\kappa: SL(V) \times H \rarr \proj(Z)$$
defined  by: $$\kappa(y,[z])=[\zeta(y)(z)]$$
is surjective onto $\proj(Z)^{\barr{e}UN}$.

The unstable locus, $\proj(Z)^{UN}$, is closed. Also,
$\proj(Z)^{UN}$ is the countable union of the $\proj(Z)^{\barr{e}UN}$. Over
an uncountable algebraically closed field, any algebraic variety
that is countable union of constructible subsets is
actually the union of finitely many of them.  Therefore a
finite set of $a$-tuples, $\cal{P}$, with the
demanded property exists in the case $\bold{C}$ is uncountable.

There always exists a field extension $\bold{C} \rarr \bold{C}'$ where
$\bold{C}'$ is an uncountable algebraically closed field. By base
extension, $$\zeta_{\bold{C}'}:
SL_{\bold{C}'}(V\otimes_{\bold{C}} \bold{C}') \rarr
SL_{\bold{C}'}(Z\otimes_{\bold{C}} \bold{C}').$$
Since $\bold{C}$ is algebraically closed, it is easy to see that
the $\bold{C}$-valued closed points of
$\proj_{\bold{C}'}(Z\otimes_{\bold{C}}
\bold{C}')^{\barr{e}UN}$ are simply the points of $\proj_{\bold{C}}
(Z)^{\barr{e}UN}$. Hence, the assertion for $\bold{C}'$ implies the
assertion for $\bold{C}$. This settles the general case.
\end{pf}

\begin{pf}[Of Lemma (\ref{jill})] Let $\cal{P}$
be determined by Lemma (\ref{tom})
for the representation $\zeta$.  Let
$M_\omega$ be determined by Lemma (\ref{wally}) for the representation
$\omega$. Let $N_\zeta$ satisfy
$$ \forall \barr{e}\in \cal{P}, \ \ \ N_\zeta > |\barr{e}|.$$
Let $k_{SS}=M_\omega \cdot N_\zeta$. Suppose
$k>k_{SS}$. For each element
$$[z]\times[w] \in \proj(Z)^{UN} \times \proj(W),$$
we must show that $[z^k\otimes w]$ is unstable for the
$Sym^k(\zeta)\otimes \omega$-induced linearization on
$\proj(Sym^k(Z) \otimes W)$. Since $[z]\in \proj(Z)^{UN}$,
there exists an $\barr{e}\in \cal{P}$ such that $[z]$ is $\barr{e}$-unstable
for the $\zeta$-induced linearization on $\proj(Z)$.
Let $\lambda: \bold{C}^* \rarr SL(V)$  be a $1$-parameter
subgroup with diagonalized exponents $\barr{e}$ that
destabilizes $[z]$. Let $(\barr{z}, \barr{f})$ and $(\barr{w},
\barr{h})$ be diagonalizing data for $\zeta \circ \lambda$ and
$\omega \circ \lambda$. Let $z=\sum_{1}^{b}\xi_s\cdot  z_s$ and
$w=\sum_{1}^{c}\sigma_t\cdot  w_t$ be the basis expansions. Since
$\lambda$ destabilizes $[z]$, we see
\begin{equation}
\label{hah}
\xi_s \neq 0 \ \ \Rightarrow \ \ f_s>0.
\end{equation}
A diagonalizing basis of $Sym^k(\zeta) \circ \lambda$ can be
constructed by taking homogeneous monomials of degree $k$ in
$\barr{z}$. Denote this basis with the corresponding exponents  by
$(\barr{z}^*, \barr{f}^*)$.  Then $\barr{z}^*\otimes \barr{w}$
is a diagonalizing basis of
$$(Sym^k(\zeta)\otimes \omega)  \circ \lambda.$$
We must show that every nonzero coefficient of the
expansion of $z^k\otimes w$ in the basis
$\barr{z}^*\otimes \barr{w}$ corresponds to a positive
exponent. Suppose the basis element $z_{s^*}^*\otimes w_t$ has
a nonzero coefficient. The element $z_{s^*}^*$ must
correspond to a homogeneous polynomial of degree $k$ in
those $z_s$ for which $\xi_s\neq0$. Therefore, by (\ref{hah}), the
exponent $f_{s^*}^*$ is not less  than $k$.
The exponent corresponding to $z_{s^*}^*\otimes w_t$ is
$f_{s^*}^* + h_t$. Since
$$ |h_t| \leq |\barr{h}| < M_\omega \cdot |\barr{e}_p| <
M_\omega \cdot N_\zeta = k,$$
$$f_{s^*}^* + h_t>0.$$
The proof is complete.
\end{pf}

\section {The Construction of $\barr{U_g(e,r)}$}
\label{conlo}
\subsection{Uniform Rank}
Define $$Q_g^r(\mu,n,f_{e,r}) \subset Q_g(\mu,n,f_{e,r})$$
to be the subset corresponding to quotients
$$\tp \oh_C  \rarr E  \rarr 0$$
where $E$ has uniform rank $r$ on $C$.
Certainly $Q_g^r(\mu,n,f_{e,r})$ is $SL_{N+1}\times SL_n$ -invariant.
\begin{lm}
$Q_g^r(\mu,n,f_{e,r})$ is open and
closed in $Q_g(\mu,n, f_{e,r})$. ($Q_g^r(\mu,n,f_{e,r})$ is a union of
connected components).
\end{lm}
\begin{pf}
Let $\kappa: \cal{C} \rarr \cal{B}$ be a projective,
 flat family of Deligne-Mumford
stable genus $g$ curves over an irreducible curve. Let
$\cal E$ be  a $\kappa$-flat coherent sheaf of constant Hilbert polynomial
$f_{e,r}$ (with respect to $\om_{\cal{C}/\cal{B}}^{10})$. Suppose
there exists a $b^*\in \cal{B}$ such that $\cal{E}_{b^*}$ has uniform rank
$r$ on $\cal{C}_{b^*}=C$. Let $\{ \cal{C}_i \}$ be the irreducible
components of $\cal{C}$. Since $\kappa:\cal{C}_i \rarr \cal{B}$
is surjective of relative dimension $1$,
 each $\cal{C}_i$ contains a component of $C$. Since
the function $r(z)=dim_{k(z)}(\cal{E}\otimes k(z))$ is  upper
semicontinuous
on $\cal{C}_i$,
there is an open set $U_i\subset \cal{C}_i$ where
$r(z) \leq r$. It follows there exists an
open set $U\subset \cal{B}$
such that $\forall b\in U$, the rank of $\cal{E}_b$ on each
component of  $\cal{C}_b$ is
at most $r$. If $\exists b'$ such that $\cal{E}_{b'}$ is not of uniform
rank $r$, then (by semicontinuity)  $\exists i $ so that $r(z)$ is strictly
less than $r$ on an open $W\subset\cal{C}_i$. For $b$ in the
nonempty intersection $U \cap \kappa(W)$, ranks of $\cal{E}_b$ are
at most $r$ on each component and strictly less
than $r$ on at least one component. By equation (\ref{oiler}) of
section (\ref{linto}),
$\cal{E}_b$ can not have Hilbert polynomial $f_{e,r}$. Thus
$\forall b \in \cal{B}$, $\cal{E}_b$ has uniform rank $r$.
The Lemma is proven.
\end{pf}

\subsection{Determination of the Semistable Locus}
Select $e>e(g,r)$ and $t>t(g,r,e)$. As usual, let $n=f_{e,r}(0)$.
Let $$Z=\bigwedge^{h(\barr{s})}H^0(\proj^N,\oh_{\proj^N}(\barr{s}))^* ,$$
$$W= \bigwedge^{f_{e,r}(t)}(\tp H^0(\proj^N,\oh_{\proj^N}(t)))^*.$$
Consider the immersion
$$j_{\barr{s},t}: Q_g^r(\mu,n,f_{e,r})
 \rarr \proj(\bigwedge^{h(\barr{s})}H^0(\proj^N,\oh_{\proj^N}(\barr{s}))^*)
\times \proj(\bigwedge^{f_{e,r}(t)}(\tp H^0(\proj^N,\oh_{\proj^N}(t)
))^*)$$
defined in section (\ref{said}). Recall $\barr{s}$ is the linearization
determined by Gieseker. There are three group actions to
examine. In what follows, the superscripts $\{ S',SS' \}$ will
denote stability and semistability with respect to the
$SL_{N+1}$-action. Similarly,  $\{S'',SS''\}$ will correspond to the
$SL_n$-action,  and $\{S,SS\}$ will correspond to the $SL_{N+1} \times
SL_n$-action.

The strategy for obtaining the desired $SL_{N+1}\times SL_n$-semistable
locus is as follows.
Consider first the $SL_{N+1}$-action. For
suitable linearization, it will be shown that
$Q_g^r(\mu,n,f_{e,r})$ is contained in the
$SL_{N+1}$-stable locus and is closed in $SL_{N+1}$-semistable locus.
This assertion is a consequence of Gieseker's conditions on $H_g$ (
(i), (ii) of section (\ref{said}) )
and the results of section (\ref{abgit}). Next,
$\barr{U_g(e,r)}$ defined as the G.I.T. quotient of
$Q_g^r(\mu,n,f_{e,r})^{SS}$ by $SL_{N+1}\times SL_{n}$.
$\barr{U_g(e,r)}$ is a projective variety. Finally,
in Proposition (\ref{lizzy}), it is shown that the $SL_n$ and
$SL_{N+1}\times SL_n$ semistable loci coincide on $Q_g^r(\mu,n,f_{e,r})$.
Similarly, the stable loci coincide. The results on the
fiberwise G.I.T. problem now yield a geometric identification
of the stable and semistable loci for the $SL_{N+1}\times SL_n$-
G.I.T. problem.

By Propositions (\ref{jack}) and (\ref{jill}),
an integer $k>\{ k_{S'},k_{SS'} \}$
can be found so that:
\begin{equation}
\label{one}
\rho_Z^{-1}(\proj(Z)^{S'}) \subset (\proj(Z) \times \proj(W))_{[k,1]}^{S'},
\end{equation}
\begin{equation}
\label{two}
(\proj(Z) \times \proj(W))_{[k,1]}^{SS'} \subset
\rho_Z^{-1}(\proj(Z)^{SS'}).
\end{equation}
By (i) of section (\ref{said}), $H_g \subset \proj(Z)^{S'}$. Now
(\ref{one}) above yields
\begin{equation}
\label{skup}
H_g \times \proj(W) \subset
 (\proj(Z) \times \proj(W))_{[k,1]}^{S'}.
\end{equation}
 Therefore,
\begin{equation}
\label{three}
Q_g^r(\mu,n,f_{e,r}) \subset (\proj(Z) \times \proj(W))_{[k,1]}^{S'}.
\end{equation}
By (ii) of section (\ref{said}), $H_g$ is closed in $\proj(Z)^{SS'}$. Hence
$H_g\times \proj(W)$ is closed in $\rho_Z^{-1}(\proj(Z)^{SS'})$.
By
(\ref{two}) and the projectivity of
$Q_g^r(\mu,n,f_{e,r})$ over $H_g$:
$$Q_g^r(\mu,n,f_{e,r})  \mbox{  is closed in  }
 (\proj(Z) \times \proj(W))_{[k,1]}^{SS'}.$$
Since
$$(\proj(Z) \times \proj(W))_{[k,1]}^{SS} \subset
(\proj(Z) \times \proj(W))_{[k,1]}^{SS'},$$
it follows that
\begin{equation}
\label{four}
Q_g^r(\mu,n,f_{e,r})_{[k,1]}^{SS} \mbox{  is closed in  }
(\proj(Z) \times \proj(W))_{[k,1]}^{SS}.
\end{equation}
We define $$\barr{U_g(e,r)} = Q_g^r(\mu,n,f_{e,r})_{[k,1]}^{SS} / (SL_{N+1}
\times SL_n).$$
By (\ref{four}), $\barr{U_g(e,r)}$ is a projective variety.

We now identify the locus $Q_g^r(\mu,n,f_{e,r})_{[k,1]}^{SS}$.
Certainly
$$Q_g^r(\mu,n,f_{e,r})_{[k,1]}^{S}
\subset Q_g^r(\mu,n,f_{e,r})_{[k,1]}^{S''} ,$$
$$Q_g^r(\mu,n,f_{e,r})_{[k,1]}^{SS}
\subset Q_g^r(\mu,n,f_{e,r})_{[k,1]}^{SS''} .$$
In fact:
\begin{pr}
\label{lizzy}
There are two equalities:
$$Q_g^r(\mu,n,f_{e,r})_{[k,1]}^{S}  = Q_g^r(\mu,n,f_{e,r})_{[k,1]}^{S''} ,$$
$$Q_g^r(\mu,n,f_{e,r})_{[k,1]}^{SS}  = Q_g^r(\mu,n,f_{e,r})_{[k,1]}^{SS''} .$$
\end{pr}
\begin{pf}
We apply the Numerical Criterion.
Let
$$ \zeta^k: SL_{N+1} \times SL_n \rarr SL_{N+1} \rarr SL(Sym^k(Z)).$$
$$ {\omega}: SL_{N+1}\times SL_n \rarr SL(W) $$
denote the two representations.
 Let $\xi \in
 Q_g^r(\mu,n,f_{e,r})_{[k,1]}^{S''}$. Recalling the morphisms defined
in section (\ref{van}),
$$j_{\barr{s},t}(\xi)= (\pi(\xi), i_t(\xi)).$$
Let $\lambda:\bold{C}^* \rarr SL_{N+1}\times SL_n $ be a nontrivial
$1$-parameter subgroup given by
$$\lambda_1 : \bold{C}^* \rarr SL_{N+1},$$
$$\lambda_2: \bold{C}^* \rarr SL_n.$$
Let $\{ m_i \}$ be a diagonalizing basis for
${\zeta}^k \circ \lambda$ with weights $\{ w(m_i) \}$. Let
$\{ n_j \}$ be a diagonalizing basis for
the $\bold{C}^* \times \bold{C}^*$ representation
${\omega }\circ (\lambda_1 \times \lambda_2)$.
 Let $w_1(n_j)$ and $w_2(n_j)$ denote
the weights of the induced $\bold{C}^*$ representations
${\omega} \circ (\lambda_1 \times 1)$ and ${\omega} \circ
(1 \times \lambda_2)$.
The weights $\{w(n_j) \}$ of the $\bold{C}^*$ representation
${\omega}\circ \lambda$ are given by
$w(n_j)=w_1(n_j)+w_2(n_j).$ Finally let
$\{\barr{m}_i \}$ and $\{ \barr{n}_j \}$
denote the elements of the diagonalizing bases that appear
with nonzero coefficient in the expansions of
$\pi(\xi)$ and $i_t(\xi)$. There are three cases.
\begin{enumerate}
\item $\lambda_1=1$. Since $\xi$ is a stable point for the
$SL_n$-action, there is a $\barr{n}_j$ with $w_2(\barr{n}_j)<0.$
We see $$ w(\barr{m}_i \otimes \barr{n}_j)=w(\barr{m}_i) +
w_1(\barr{n}_j) + w_2(\barr{n}_j) =w_2(\barr{n}_j) <0 $$
for any $\barr{m}_i$.
\item $\lambda_2=1$. By (\ref{three}), $\xi$ is a stable point for the
$SL_{N+1}$-action. Hence there exists a pair $\barr{m}_i$, $\barr{n}_j$
so that
$$w(\barr{m}_i \otimes \barr{n}_j) =w(\barr{m}_i)+w_1(\barr{n}_j) <0.$$
\item $\lambda_1 \neq 1, \ \lambda_2 \neq 1$. Since $\xi$ is a stable
point for the
$SL_n$-action, there is a $\barr{n}_j$ with $w_2(\barr{n}_j)<0.$
By (\ref{skup}), $(\pi(\xi) \otimes \barr{n}_j)$ is a stable point
for the $SL_{N+1}$-action. Hence there exists an element  $\barr{m}_i$
so that
$$w(\barr{m}_i)+w_1(\barr{n}_j) <0.$$
Therefore,
$$ w(\barr{m}_i \otimes \barr{n}_j)=w(\barr{m}_i) +
w_1(\barr{n}_j) + w_2(\barr{n}_j) <0.$$
\end{enumerate}
By the Numerical Criterion, $\xi\in Q_g^r(\mu,n,f_{e,r})_{[k,1]}^S.$
The proof for the semistable case is identical.
\end{pf}
By Theorem (\ref{fred}), we see the points of
$Q_g^r(\mu,n,f_{e,r})_{[k,1]}^{SS}$ correspond exactly to
quotients
$$\tp \oh_C \rarr E \rarr 0$$
where $E$ is slope-semistable, torsion free sheaf of uniform rank
$r$ on a $10$-canonical, Deligne-Mumford stable, genus $g$ curve
 $C\subset\proj^N$ and
$$\psi: \tp H^0(C,\oh_C) \rarr H^0(C,E)$$
is an isomorphism.
Similarly for $Q_g^r(\mu,n,f_{e,r})_{[k,1]}^S$.

We now examine orbit closures. Suppose
 $\barr{\xi} \in Q_g^r(\mu,n,f_{e,r})_{[k,1]}^{SS}$ lies in the
$SL_{N+1} \times SL_n$-orbit closure
of of  $\xi  \in Q_g^r(\mu,n,f_{e,r})_{[k,1]}^{SS}$.
Let $$\gamma=(\gamma_1,\gamma_2): \bigtriangleup
-\{p\} \rarr SL_{N+1} \times SL_n$$
 be a morphism of a nonsingular, pointed curve such that
$$Lim_{z\rarr p}(\gamma(z)\cdot\xi) = \barr{\xi}.$$
 It follows that
$$Lim_{z\rarr p}(\gamma_1(z)\cdot \pi(\xi)) = \pi(\barr{\xi}).$$
Since $H_g \subset \proj(Z)^{S}$,  $\pi(\barr{\xi})$
lies in the $SL_{N+1}$-orbit of $ \pi(\xi)$. After a possible base
change, we can assume $\gamma_1$ extends over $p\in \bigtriangleup$ to
$$\barr{\gamma_1}: \bigtriangleup_p \rarr SL_{N+1}.$$
Since $Q_g^r(U_{\pi(\xi)}, n,f_{e,r})$ is projective,
$$\mu : \bigtriangleup - \{p \} \rarr Q_g^r(U_{\pi(\xi)}, n,f_{e,r})$$
defined by $\mu(z)= \gamma_2(z) \cdot \xi$ extends to
$\bigtriangleup_p$. Let
$\hat{\xi}= Lim_{z\rarr p}(\gamma_2(z)\cdot \xi)= \mu(p).$
By considering the map
 $$\barr{\gamma_1}\cdot \mu:\bigtriangleup_p \rarr Q_g^r(\mu,n,f_{e,r})$$
defined by
$$\barr{\gamma_1}\cdot \mu(z)= \barr{\gamma_1}(z)\cdot \mu(z) $$
we obtain,
$$\barr{\gamma_1}(p) \cdot \hat{\xi}= \barr{\gamma_1}\cdot \mu (p)
 = Lim_{z\rarr p}(\gamma_1(z)\cdot \gamma_2(z)\cdot \xi) =$$
$$Lim_{z\rarr p}(\gamma(z)\cdot \xi) =\barr{\xi}.$$
We have shown the $SL_{N+1}$-orbit of $\barr{\xi}$  meets
the $SL_n$-orbit closure of $\xi$. If $\xi$, $\barr{\xi}$ corresponds to
a slope-semistable, torsion free quotients $E$,  $\barr{E}$ on
$C$, $\barr{C} \subset \proj^N$, certainly $C$, $\barr{C}$ must be
projectively equivalent. The elements of the $SL_{N+1}$-orbit of
$\barr{\xi}$ that lie over $C$ are simply the images of $\barr{E}$
under automorphisms of $C$. Now from section (\ref{sesh}), we
conclude two semistable orbits $\xi$ and $\barr{\xi}$ are
identified in the quotient $\barr{U_g(e,r)}$ if and only if
$$\pi(\xi) \equiv \pi(\barr{\xi}) \equiv  [C]$$ and
the corresponding semistable, torsion free quotient
sheaves $E$, $\barr{E}$ have Jordan-Holder factors that
differ by an automorphism of $C$. We see:
\begin{tm}
\label{priya}
$\barr{U_g(e,r)}$ parametrizes aut-equivalence classes of
slope-semistable, torsion free sheaves of uniform rank $r$
and degree $e$ on Deligne-Mumford stable curves of genus $g$.
\end{tm}
\noindent
Finally, since $$ Q_g^r(\mu,n,f_{e,r})_{[k,1]}^{SS} \rarr H_g \rarr \mgb$$
is an $SL_{N+1}\times SL_{n}-$ invariant morphism, there exists
a map
$$\eta: \barr{U_g(e,r)} \rarr \mgb.$$

\section{Basic Properties of $\barr{U_g(e,r)}$}
\label{laslo}
\subsection{The Functor}
Let $\barr{\cal{U}_g(e,r)}$ be the functor that associates to
each scheme $S$ the set of equivalence classes of the following
data:
\begin{enumerate}
\item A flat family of Deligne-Mumford stable genus $g$ curves,
$\mu: \cal{C} \rarr S$.
\item A $\mu$-flat coherent sheaf $\cal{E}$ on $\cal{C}$ such that:
\begin{enumerate}
\item [(i.)] $\cal{E}$
is of constant Hilbert polynomial $f_{e,r}$ (with respect to
$\om_{\cal{C}/S}^{10}$).
\item [(ii.)] $\cal{E}$ is a slope-semistable, torsion free
sheaf of uniform rank $r$ on each fiber.
\end{enumerate}
\end{enumerate}
Two such data sets are equivalent if there exists an $S$ isomorphism
$\phi: \cal{C} \rarr \cal{C}'$
and a line bundle $L$ on $S$ so that
$\cal{E} \cong  \phi^*(\cal{E}') \otimes \mu^*L$.
\begin{tm}
\label{thwee}
There exists  a natural transformation
$$\phi_U: \barr{\cal{U}_g(e,r)} \rarr Hom(*, \barr{U_g(e,r)}).$$
$\barr{U_g(e,r)}$ is universal in the following sense.
If $Z$ is a scheme and
$\phi_Z: \barr{\cal{U}_g(e,r)} \rarr Hom( *, Z)$ is
a natural transformation, there exists a unique morphism
$ \gamma: \barr{{U}_g(e,r)} \rarr Z $ such that
the transformations $\phi_Z$ and $\gamma\circ \phi_U$ are
equal.
\end{tm}

\begin{pf}
Let $e>e(g,r)$.
Let $\mu: \cal{C} \rarr S$ and $\cal{E}$ on $C$ satisfy (1) and (2)
above. Note $\mu_*(\om_{\cal{C}/S}^{10})$ is a locally free
sheaf of rank $N+1= 10(2g-2)-g+1$. Since $\cal{E}_s$ is
nonspecial for each $s\in S$,
$\mu_*(\cal{E})$ is locally
free of rank $n=f_{e,r}(0)$. Choose an open cover
$\{W_i \}$ of $S$ trivializing  both $\mu_*(\om_{\cal{C}/S}^{10})$ and
$\mu_*(\cal{E})$.
Let $V_i=\mu^{-1}(W_i)$. For each $i$, we obtain isomorphisms:
\begin{equation}
\label{hoot}
\bold{C}^{N+1}\otimes \oh_{W_i} \cong \mu_*(\om_{\cal{C}/S}^{10})|_{W_i}
\end{equation}
\begin{equation}
\label{toot}
 \bold{C}^{n}\otimes \oh_{W_i} \cong \mu_*(\cal{E})|_{W_i}.
\end{equation}
These isomorphisms yield surjections:
$$ \bold{C}^{N+1}\otimes \oh_{V_i} \cong
 \mu^*\mu_*(\om_{\cal{C}/S}^{10})|_{V_i} \rarr \om_{\cal{C}/S}^{10}
|_{V_i} \rarr 0$$
$$ \tp \oh_{V_i} \cong
 \mu^*\mu_*(\cal{E})|_{V_i} \rarr \cal{E}
|_{V_i} \rarr 0$$
The first surjection embeds $V_i$ in $W_i \times \proj^N$.
By the universal property of the Quot scheme $Q_g(\mu,n,f_{e,r})$
and the second surjection,
there exists a map
$$\phi_i: W_i \rarr Q_g(\mu,n,f_{e,r}).$$
For $t>t(g,r,e)$, $\phi_i(W_i) \subset Q_g^r(\mu,n,f_{e,r})_{[k,1]}^{SS}$.
 On the overlaps $W_i\cap W_j$, $\phi_i$ and
$\phi_j$ differ by a morphism
$$W_i \cap W_j \rarr PSL_{N+1} \times PSL_n$$ corresponding
to the choice  of trivialization in (\ref{hoot}) and (\ref{toot}).
Hence there exists a well defined morphism
$$\phi: S \rarr \barr{U_g(e,r)}.$$
The functoriality of the universal property of the Quot scheme implies
$\phi$ is functorially associated to
the data $\cal{E}$ and $\mu:\cal{C} \rarr S$.
We have shown there exists a natural transformation
$$\phi_U:\barr{\cal{U}_g(e,r)} \rarr Hom( *, \barr{U_g(e,r)}).$$

Suppose
$\phi_Z:\barr{\cal{U}_g(e,r)} \rarr Hom( *, Z)$ is a
natural transformation. There exists a canonical
element of
$\delta \in \barr{\cal{U}_g(e,r)}(Q_g^r(\mu,n,f_{e,r})_{[k,1]}^{SS})$
corresponding
to the universal family on the Quot scheme. The morphism
$$\phi_Z(\delta) :Q_g^r(\mu,n,f_{e,r})_{[k,1]}^{SS} \rarr Z$$
is $SL_{N+1} \times SL_{n}$-invariant. Hence $\phi_Z(\delta)$
descends to $\gamma: \barr{U_g(e,r)}\rarr Z$.
The two transformations $\phi_Z$ and $ \gamma \circ \phi_U$
agree on $\delta$. Naturality now implies $\phi_Z= \gamma \circ \phi_U$.
\end{pf}

By previous considerations, there are natural transformations
$$t_l:\barr{\cal{U}_g(e,r)} \rarr \barr{\cal{U}_g(e+rl(2g-2),r)}$$
given by $\cal{E} \rarr \cal{E}\otimes \om_{\cal{C}/S}^l$.
By Theorem (\ref{thwee}), these induce  natural isomorphisms
$$t_l:\barr{U_g(e,r)} \rarr \barr{U_g(e+rl(2g-2),r)}.$$
The arguments in the above proof imply a useful Lemma:

\begin{lm}
\label{openy}
Let $\mu: \cal{C} \rarr S$ be a flat family of Deligne-Mumford
stable, genus $g\geq2$ curves. Let $\cal{E}$ be a $\mu$-flat
coherent sheaf
on $\cal{C}$. The condition that $\cal{E}_s $ is
a slope-semistable torsion free  sheaf of uniform rank on $\cal{C}_s$ is
open on $S$.
\end{lm}
\begin{pf}
Suppose $\cal{E}_{{s_0}} $ is
a slope-semistable sheaf of uniform rank $r$ on $\cal{C}_{{s}_0}$
for some ${s}_0\in S$.
There exists a integer $m$ so :
\begin{enumerate}
\item [(i.)] $h^1(\cal{E}_s\otimes \om_{\cal{C}_s}^m, \cal{C}_s)=0$ for all
          $s\in S$.
\item [(ii.)] $\cal{E}_s\otimes \om_{\cal{C}_s}^m$ is generated by global
section for all $s\in S$.
\item [(iii.)] $degree(\cal{E}_{s_0}\otimes \om_{\cal{C}_{s_0}}^m) > e(g,r)$.
\end{enumerate}
It suffices to prove the Lemma
for $\cal{E} \otimes \om_{\cal{C}/S}^m$. Let $f_{e,r}$ be
the Hilbert polynomial of $\cal{E} \otimes \om_{\cal{C}/S}^m$.
By the proof of Theorem (\ref{thwee}), there
exists an open set $W \subset S$ containing $s_0$  and
a morphism $$\phi: W \rarr Q_g(\mu,n=f_{e,r}(0), f_{e,r})$$
such that $\cal{E}\otimes \om_{\cal{C}/S}^m$ is isomorphic to
the $\phi$-pull back of the universal quotient.
Since $\phi({s_0}) \in Q_g^r(\mu,n,f_{e,r})_{[k,1]}^{SS}$ and
the latter is open, the Lemma is proven.
\end{pf}

\subsection{Deformations of Torsion Free Sheaves
and the Irreducibility of $\barr{U_g(e,r)}$}
We study deformation properties of uniform rank,
torsion free sheaves on nodal curves.
\begin{lm}
\label{d1}
Let
$\mu: S \rarr Spec(\bold{C}[t])$ be a $\mu$-flat, nonsingular
surface.
If $\cal{E}$ is a $\mu$-flat sheaf on $S$ such that the
restriction $\cal{E}/t\cal{E} =\cal{E}_0$
is torsion free, then $\cal{E}_0$ is locally free.
\end{lm}
\begin{pf}
Let $z\in \mu^{-1}(0)$. Since $S$ is a nonsingular surface, the
local ring $\oh_{S,z}$ is regular of dimension $2$.
Consider the $\oh_{S,z}$-module $\cal{E}_z$. Since $\cal{E}$
is $\mu$-flat, $t$ is a $\cal{E}_z$-regular element.
Since $\cal{E}_0=\cal{E}/t\cal{E}$ is torsion free,
$depth_{\oh_{S,z}}(\cal{E}_z)\geq 2$. We have the
Auslander-Buchsbaum relation ([M]):
$$proj.\  dim_{\oh_{S,z}}(\cal{E}_z) + depth_{\oh_{S,z}}(\cal{E}_z) =
dim(\oh_{S,z})=2.$$
We conclude $\ proj.\ dim_{\oh_{S,z}}(\cal{E}_z)=0$.
Hence $\cal{E}_z$ is free over $\oh_{S,z}$. It follows that
$\cal{E}_0$ is locally free.
\end{pf}
Lemma (\ref{d1}) shows that it is not possible to deform a
torsion free, non-locally free sheaf on a nodal curve
to a locally free sheaf on a nonsingular  curve if the deformations
at the nodes have local equations of the form
$(xy-t)\subset Spec(\bold{C}[x,y,t])$.
However, the next Lemma shows such deformations exist locally if the
deformations of the nodes have local equations of the form
$(xy-t^2)$. In Lemma (\ref{deform}), it is shown these local deformations
can be globalized.

\begin{lm}
\label{d2}
Let $S\subset Spec(\bold{C}[x,y,t])$ be the subscheme defined
by the ideal $(xy-t^2)$. Let
$\mu: S \rarr Spec(\bold{C}[t])$. Let $\zeta=(0,0,0)\in S$. There exists a
$\mu$-flat sheaf $\cal{E}$ on $S$ such that
$\cal{E}_{t\neq 0}$ is locally free and
$\cal{E}_0 \cong m_{\zeta}$ where $m_{\zeta}$ is the
maximal ideal defining $\zeta$ on $S_0$.
\end{lm}

\begin{pf}
There exists a section $L$ of $\mu$ defined by
the ideal $(x-t,y-t)$. Let $\cal{E}$ be the coherent
sheaf corresponding to this ideal.  We have the
exact sequence:
\begin{equation}
\label{pddy}
0\rarr\cal{E} \rarr \oh_S \rarr \oh_L \rarr 0.
\end{equation}
Since $\oh_S$ is torsion free over $\bold{C}[t]$,
so is $\cal{E}$. $\cal{E}$ is therefore $\mu$-flat.
Since $\oh_L$ is $\mu$-flat, sequence (\ref{pddy})
remains exact after restriction to the special fiber. Hence
$\cal{E}_0\cong m_\zeta$.
\end{pf}

\begin{lm}
\label{deform}
Let $C$ be a Deligne-Mumford stable curve of genus $g\geq2$.
Let $E$ be a slope-semistable torsion free sheaf of uniform
rank $r$ on $C$.  Then there exists a family
$\mu: \cal{C} \rarr \bigtriangleup_0$ and a $\mu$-flat coherent sheaf
$\cal{E}$ on $\cal{C}$ such that :
\begin{enumerate}
\item $\bigtriangleup_0$ is a pointed curve.
\item $\cal{C}_0 \cong C$,  $\ \ \forall t\neq 0$ $\ \ \cal{C}_t$ is
a complete, nonsingular, irreducible genus $g$ curve.
\item $\cal{E}_0 \cong E$,  $\ \ \forall t\neq 0$ $\ \ \cal{E}_t$ is
a slope-semistable torsion free sheaf of rank $r$.
\end{enumerate}
\end{lm}

\begin{pf}
Let $z\in C$ be a node.
Since $E$ is torsion free and of uniform rank $r$, it
follows from Propositions (2) and (3) of chapter (8) of
[Se]:
\begin{equation}
\label{locy}
E_z \cong \bigoplus_{1}^{a_z} \oh_{C_z} \oplus \bigoplus _{a_z+1}^{r}
m_z
\end{equation}
where $m_z$ is the localization of the ideal of the node $z$.
To simplify the deformation arguments, let
$\bold{C}$ be the field of complex numbers.
Let $B(d)\subset \bold{C}^2$ be the open ball of radius $d$
with respect to the Euclidean norm;
let $B(d_1,d_2)\subset \bold{C}^2$ be the open annulus.
Disjoint open Euclidean neighborhoods, $z\in U_z\subset C$, can be
chosen for each node of $C$ satisfying:
\begin{enumerate}
\item[(i.)] $U_z$ is analytically
isomorphic to
$B(d_z) \cap (xy=0)\subset \bold{C}^2$.
\item[(ii.)]
$E|_{U_z}\cong \bigoplus_{1}^{a_z} \oh_{U_z} \oplus \bigoplus _{a_z+1}^{r}
m_z$.
\end{enumerate}
Let $V_z\subset U_z$ be the closed neighborhood of radius $d_z/2$.
Let $W = C \setminus \cup V_z$.
A deformation of $C$ can be given the by the open cover:
$$\{W\times \bigtriangleup_0\}\cup \{Def_z |z \in C_{ns}\}$$
where $Def_z$ (to be defined below) is an open subset of
$$B(d_z)\times\bigtriangleup_0\  \cap\  (xy-t^2=0)\ \ \ \
\subset \ \ \bold{C}^2\times
\bigtriangleup_0$$
containing $(0,0,0)$. $Def_z$ is a local smoothing at the node $z$.
Define
$$K_z \ \ = \  \ B\paren{{d_z\over 2}, {d_z\over 2}+\epsilon_z}
\times
\bigtriangleup_0 \ \cap\  (xy=0)\ \  \subset \ \ \bold{C}^2\times
\bigtriangleup_0.$$
Note $B(r,s)$ is the
annulus.
Let $\mu$ denote the projection to $\bigtriangleup_0$.
For $\epsilon_z>0$ (small with respect to $\delta_z$)
and $|t|< \delta_z$, it is not hard
to find an isomorphism $\gamma_z$ commuting with $\mu$:
$$\gamma_z: K_z
\rarr D_z \ \ \subset B(d_z)\times\bigtriangleup_0\  \cap\  (xy-t^2=0)$$
such that
\begin{equation}
\label{duder}
B\paren{{d_z\over 3}}\times \bigtriangleup_0 \ \cap \ D_z = \emptyset
\end{equation}
and
$\gamma_z$  extends the identity on $t=0$.
Such a $\gamma_z$ can be constructed by considering
the holomorphic flow of a algebraic vector field on $(xy-t^2=0)$.
The space
$$B(d_z)\times\bigtriangleup_0 \ \cap\  (xy-t^2=0)\ \  \setminus \ D_z $$
is disconnected. Let $Def_z$ be the complement of the
component not containing $(0,0,0)$.
 The isomorphism $\gamma_z$
determines a patching of $W\times \bigtriangleup_0$ and $Def_z$
along $K_z \simeq D_z$ in
the obvious manner. The constructed family satisfies claims
(1) and (2) of the Lemma.

$E_0$ can be extended trivially on $W\times \bigtriangleup_0$ to
$\cal{E}|_W$.
$\cal{E}|_{K_z}$ is trivial by condition (ii).
By Lemma (\ref{d2}), $m_z$ can be flatly  extended to a line bundle $L_z$
on $Def_z$. By (\ref{duder}) and the construction of Lemma (\ref{d2}),
$L_z$ can be assumed to be
 trivial on $D_z$. By patching
$$\bigoplus_{1}^{a_z} \oh_{Def_z} \oplus \bigoplus _{a_z+1}^{r}
L_z$$
 along $K_z\simeq D_z$,
$\cal{E}$ can be defined such that
$\cal{E}_0\cong E$ and $\cal{E}_{t\neq 0}$ is locally free.
Indeed, such a patching exists for $t=0$ by condition (ii).
The patching can be extended trivially along $K_z$ since
$$K_z= K_{z,t=0} \times \bigtriangleup_0.$$
Now condition (3) follows by Lemma (\ref{openy}).
For a general ground field, the
\'etale topology must be used.
\end{pf}

\begin{pr}
$\barr{U_g(e,r)}$ is an irreducible variety.
\end{pr}
\begin{pf}
Consider the morphism
$\pi_{SS}: Q^r_g(\mu,n=f_{e,r}(0),f_{e,r})_{[k,1]}^{SS} \rarr H_g$.
By Proposition (24) of chapter (1) of [Se], the
scheme
 $$\pi_{SS}^{-1}([C])=Q_g(C,n, f_{e,r})^{SS}$$
is irreducible for each nonsingular curve $C$,
 $[C]\in H_g$. Since the locus $H_g^0\subset H_g$ of nonsingular
curves is irreducible, $\pi_{SS}^{-1}(H_g^0)$ is irreducible.
By Lemma (\ref{deform}), $\pi_{SS}^{-1}(H_g^0)$ is dense in
$Q^r_g(\mu,n,f_{e,r})_{[k,1]}^{SS}$.
 Finally, since there is a surjective morphism
$$Q^r_g(\mu,n,f_{e,r})_{[k,1]}^{SS} \rarr \barr{U_g(e,r)},$$
we conclude $\barr{U_g(e,r)}$ is irreducible.
\end{pf}

\subsection{The Normality of $\barr{U_g(e,1)}$}
We need an infinitesimal analogue of Lemma (\ref{d1}).
First, we establish some notation. Let
 $$R= \bold{C}[[x,y]][\epsilon]/(xy-\epsilon,\epsilon^2).$$
Let $A= R/\epsilon R \cong \bold{C}[[x,y]]/(xy)$. Let
$m=(x,y)$ be the maximal ideal of $A$.
There is a natural, flat  inclusion of rings
$\bold{C}[\epsilon]/(\epsilon^2) \rarr R$. For a
$\bold{C}[\epsilon]/(\epsilon^2)$-module $M$, let $*\epsilon$
denote the $M$ endomorphism given by multiplication by $\epsilon$.
\begin{lm}
\label{ey}
There does not exist a $\bold{C}[\epsilon]/(\epsilon^2)$-flat
$R$-module $E$ such that $$E/\epsilon E \cong m$$ as  $A$-modules.
\end{lm}

\begin{pf}
Suppose such an $E$ exists. Let
$$\alpha: E \rarr E/\epsilon E \stackrel{\sim}{\rarr} m.$$
Let $e_x$, $e_y \in E$ satisfy $\alpha(e_x)=x$ and
$\alpha(e_y)=y$. We obtain a morphism
$$\beta :  R\oplus R \rarr E $$
defined by $\beta(1,0)=e_x$ and $\beta(0,1)=e_y$.
Since $\epsilon$ is nilpotent,
$\beta$ is surjective by Nakayama's Lemma. We have an exact
sequence
$$0\rarr N \rarr R \oplus R \rarr E \rarr 0.$$
Since $E$
and $R\oplus R$  are  $\bold{C}[\epsilon]/(\epsilon^2)$-flat,
$N$ is $\bold{C}[\epsilon]/(\epsilon^2)$-flat.
Hence, there is an {\em exact} sequence
\begin{equation}
\label{fyttf}
0 \rarr \epsilon N \rarr N \stackrel{* \epsilon}{\rarr} N
\end{equation}
obtained by tensoring
$0 \rarr (\epsilon) \rarr \bold{C}[\epsilon]/(\epsilon ^2)
 \stackrel{* \epsilon}{\rarr} \bold{C}[\epsilon]/(\epsilon ^2)$
with $N$.
Since $\beta \bigl( (y,0) \bigr) \in \epsilon E$,  there exists
an element $n=(y+ \epsilon f(x,y), \epsilon g(x,y))$
in $N$. Consider $x n \in N$:
$$xn = (xy+ \epsilon xf(x,y), \epsilon xg(x,y))=(\epsilon(1+
xf(x,y)), \epsilon x g(x,y)).$$
Note $\epsilon x n=0$.
By the exactness of  (\ref{fyttf}), there exists an $\barr{n}\in N$
satisfying $\epsilon \barr{n}= xn$. Since
$R\oplus R$ is flat over $\bold{C}[\epsilon]/(\epsilon^2)$,
 any such $\barr{n}$ must
be of the form
$$\barr{n}= (\ 1+xf(x,y) + \epsilon \hat{f}(x,y),\
  xg(x,y)+\epsilon \hat{g}(x,y)\ ).$$
Since
$\alpha \circ \beta (\barr{n}) = x+x^2f(x,y) \neq 0$ in $m$,
$\ \barr{n}$ can not lie in $N$.  We have a contradiction. No such $E$
can exist.
\end{pf}

To prove $\barr{U_g(e,1)}$ is normal, it suffices to
show $Q_g^1(\mu, n=f_{e,1}, f_{e,1})_{[k,1]}^{SS}$ is nonsingular.
The nonsingularity is established by computing the dimension of
$Q_g^1(\mu,n,f_{e,1})_{[k,1]}^{SS}$ and
then bounding the dimension of the Zariski tangent space at each point.
The Zariski tangent spaces are controlled by a study of the
differential $d\pi_{SS}$ where $\pi_{SS}$ is the canonical map
$$\pi_{SS}: Q_g^1(\mu, n, f_{e,1})_{[k,1]}^{SS} \rarr H_g.$$
\begin{lm}
$Q_g^1(\mu, n=f_{e,1}(0), f_{e,1})_{[k,1]}^{SS}$ is nonsingular.
\end{lm}
\begin{pf}
Consider the  universal quotient sequence over
 $Q^1_g(\mu,n,f_{e,1})_{[k,1]}^{SS}$,
$$0 \rarr \F \rarr \tp \oh_{Q^{ss}\times U} \rarr \E \rarr 0.$$
 Let
$\xi \in Q^1_g(\mu,n,f_{e,1})_{[k,1]}^{SS}$
be a closed point and let $C=U_{\pi(\xi)}$.
$\xi$ corresponds to a quotient
\begin{equation}
\label{smooth}
0 \rarr \F_\xi \rarr \tp \oh_C \rarr \E_\xi \rarr 0.
\end{equation}
There is a natural identification of
the Zariski tangent space to $Q^1_g(C,n,f_{e,1})$ at $\xi$:
$$T_\xi(Q^1_g(C,n,f_{e,1})) \cong   H^0(C,\underline{Hom}(\F_\xi,\E_\xi))$$
 (see [Gr]).
If $h^1(C,\underline{Hom}(\F_\xi,\E_\xi))=0$
and the deformations of
\begin{equation}
\label{george}
\tp \oh_C \rarr \E_\xi \rarr 0
\end{equation}
are locally
unobstructed, then
$\xi$ is a nonsingular point of $Q^1_g(C,n,f_{e,1})$.
If $\E_\xi$ is locally free, the deformations of (\ref{george})
are locally unobstructed.
Sequence (\ref{smooth})
yields:
\begin{equation}
\label{exty}
0 \rarr \underline{Hom}(\E_\xi,\E_\xi) \rarr
  \underline{Hom}(\tp \oh_C,\E_\xi)
\rarr \underline{Hom}(\F_\xi,\E_\xi) \rarr \underline{Ext}^1(\E_\xi,\E_\xi)
\rarr 0.
\end{equation}
Since $\underline{Ext}^1(\E_\xi,\E_\xi)$ is torsion and
$$h^1(C,\underline{Hom}(\tp \oh_C,\E_\xi)) =
n\cdot h^1(C,\E_\xi) =0,$$
we obtain $h^1(C,\underline{Hom}(\F_\xi,\E_\xi))=0$.

{}From (\ref{locy}), at each node $z\in C$, $\E_\xi$ is
either locally free or locally isomorphic to $m_z$.
Let $s$ be the number nodes where $\E_{\xi,z}\cong m_z$.
Using (\ref{exty}), we compute
$\chi(\underline{Hom}(\F_\xi,\E_\xi))$
in terms of $s$:
$$\chi(\underline{Hom}(\F_\xi,\E_\xi))= -\chi(\underline{Hom}(\E_\xi,\E_\xi))
+\chi( \underline{Hom}(\tp \oh_C,\E_\xi))
+\chi ( \underline{Ext}^1(\E_\xi,\E_\xi)).$$
It is clear
$\chi( \underline{Hom}(\tp \oh_C,\E_\xi))=n^2.$
Let $A=\bold{C}[[x,y]]/(xy)$ and $m=(x,y) \subset A$ as above.
It is not hard to establish:
\begin{equation}
\label{paul}
Ext_A^1(m,m)= \bold{C}^2,
\end{equation}
\begin{equation}
\label{doug}
0 \rarr A \stackrel{i}{\rarr} Hom_A(m,m) \rarr \bold{C} \rarr 0
\end{equation}
where $i$ is the natural inclusion.
Since $\underline{Ext}^1(\E_\xi,\E_\xi)$ is a torsion sheaf
supported at each $z\in C$ where $\E_{\xi,z}\cong m_z$ with
stalk (\ref{paul}),
$$\chi ( \underline{Ext}^1(\E_\xi,\E_\xi))= 2s.$$
There is a natural inclusion
$$0 \rarr \oh_C \stackrel{i}{\rarr}
\underline{Hom}(\E_\xi,\E_\xi) \rarr \delta \rarr 0.$$
Since $\E_\xi$ is of rank $1$, $\delta$ is a torsion
sheaf supported at the nodes where $\E_\xi$ is not locally free.
At these nodes, $\delta$ can be determined locally by
(\ref{doug}). Hence
$$\chi(\underline{Hom}(\E_\xi,\E_\xi))= 1-g+s.$$
Summing the Euler characteristics yields:
$$h^0(C,\underline{Hom}(\F_\xi,\E_\xi))= g-1-s+ n^2 +2s= n^2-1+g+s.$$
If $C$ is a nonsingular curve,
$\E_\xi$ is locally free on $C$. The above results
show that $\xi$ is a nonsingular point of $Q^1_g(C,n,f_{e,1})$.
Thus $dim(Q^1_g(\mu,n,f_{e,1})_{[k,1]}^{SS})= n^2-1+g+ dim(H_g)$.

Let $\xi \in Q^1_g(\mu,n,f_{e,1})_{[k,1]}^{SS}$
be a closed point with the notation given above.
We examine the exact differential sequence:
$$0\rarr
T_\xi(\pi_{SS}^{-1}[C]) \rarr T_\xi(Q^1_g(\mu,n,f_{e,1})_{[k,1]}^{SS})
\stackrel {d\pi_{SS}}{\rarr} T_{[C]}(H_g).$$
Recall $H_g$ is nonsingular.
By
Lemma (\ref{ey}) and the surjection of $T_{[C]}(H_g)$ onto
the miniversal deformation space of $C$,
$$dim( im (d\pi_{SS})) \leq dim(H_g)-s$$
($s$ is the number of nodes where $\E_\xi$ is not locally free).
By previous results:
$$dim(T_\xi(\pi_{SS}^{-1}[C]))= n^2-1+g+s.$$
It follows:
\begin{equation}
\label{zsde}
dim(T_\xi(Q^1_g(\mu,n,f_{e,1})_{[k,1]}^{SS})) \leq (n^2-1+g+s)+(dim(H_g)-s)
\end{equation}
$$=  dim(Q^1_g(\mu,n,f_{e,1})_{[k,1]}^{SS}).$$
Equality must hold in (\ref{zsde}). $\xi$ is therefore
a nonsingular point of $Q^1_g(\mu,n,f_{e,1})_{[k,1]}^{SS}$.
\end{pf}
As a consequence, we obtain:
\begin{pr}
\label{popo}
$\barr{U_g(e,1)}$ is normal.
\end{pr}

\section{The Isomorphism Between $\barr{U_g(e,1)}$ and $\barr{P_{g,e}}$}
\label{lcap}
\subsection{A Review of $\barr{P_{g,e}}$}
\label{fofo}
In this section, the compactification of the
universal Picard variety of degree $e$
line bundles, $\barr{P_{g,e}}$, described in [Ca] is considered.
Let $e$ be large enough to guarantee
the existence and properties of $\barr{P_{g,e}}$ and $\barr{U_g(e,1)}$.
A natural isomorphism
$\nu : \barr{P_{g,e}} \rarr \barr{U_g(e,1)}$ will be constructed.

We follow the notation of [Ca], [Gi].
A Deligne-Mumford {\em quasi-stable}, genus $g$ curve $C$ is a
Deligne-Mumford semistable, genus $g$ curve with destabilizing
chains of length at most one.
Let $\psi:C\rarr C_s$ be the canonical
contraction to the Deligne-Mumford
stable model.  For each
complete subcurve $D\subset C$, let $D^c=\barr{C \setminus D}$. Define:
$$k_D= \# (D \cap D^c).$$
Let $\om_{C,D}$ be the degree of the canonical bundle
$\om_C$ restricted to $D$. Let $L$ be a degree $e$ line bundle
on $C$.
Denote by $L_D$ the restriction of $L$ to $D$.
Let $e_D$ be the degree of
$L_D$.  $L$ has
{\em (semi)stable multidegree} if for each complete, proper subcurve
$D\subset C$, the following holds:
\begin{equation}
\label{basinq}
e_D   - e \cdot \paren{\om_{C,D}\over 2g-2} \ \ \  (\leq) < \ \ \  k_D/2.
\end{equation}

Consider the Hilbert scheme $H_{g,e,M}$ of degree $e$, genus $g$
curves in $\proj^M$ where $M=e-g+1$. In [Gi], it is shown there
exists an open locus
$Z_{g,e}\subset H_{g,e,M}$
parametrizing nondegenerate, Deligne-Mumford quasi-stable, genus
$g$ curves $C\subset \proj^M$ satisfying:
\begin{enumerate}
\item[(i.)] $h^1(C, \oh_C(1))= 0$.
\item[(ii.)] $\oh_C(1)$ is of semistable multidegree on $C$.
\end{enumerate} In [Ca], $\barr{P_{g,e}}$ is constructed as the
G.I.T. quotient
$$\barr{P_{g,e}} \stackrel{\sim}{=} Z_{g,e}/ SL_{M+1}$$
(for a suitable linearization). $\barr{P_{g,e}}$ is a moduli space of
line bundles of semistable multidegree on Deligne-Mumford quasi-stable curves
(up to equivalence) compactifying the universal Picard variety.

The construction of the isomorphism $\nu: \barr{P_{g,e}} \rarr
\barr{U_g(e,1)}$ proceeds as follows. If $L$ is a very ample
line bundle of semistable multidegree on a Deligne-Mumford
quasi-stable curve $C$, then $\psi_*(L)$ is a slope-semistable
torsion free sheaf of uniform rank $1$ on $C_s$. This is the
result of Lemma (\ref{ms}). The map $\nu$ is constructed by globalizing
this correspondence.
There is a universal
curve
$$U_Z \hookrightarrow Z_{g,e}\times \proj^M.$$
A deformation study shows $Z_{g,e}$ and $U_Z$ are nonsingular
quasi-projective varieties ([Ca], Lemma 2.2, p.609). There exists
a canonical contraction map $\psi:
U_Z \rarr U^s_Z$ over $Z_{g,e}$.
The map $\psi$ contracts each fiber of $U_Z$
over $Z_{g,e}$ to its Deligne-Mumford
stable model. $U^s_Z$ is a flat, projective family of Deligne-Mumford
stable curves over $Z_{g,e}$.
Let $L=\oh_{U_Z}(1)$ and
$\cal{E}= \psi_*(L)$. In Lemma (\ref{flaty}),  $\cal{E}$ is shown to be a
flat family of slope-semistable torsion free sheaves of uniform rank $1$
and degree $e$ over $Z_{g,e}$. Care is required in establishing flatness.
The argument depends upon
Zariski's theorem on formal functions and the criterion of
Lemma (\ref{tpt}). By the universal property of $\barr{U_g(e,1)}$, $\cal{E}$
induces a map $\nu_Z: Z_{g,e} \rarr \barr{U_g(e,1)}$. Since $\nu_Z$ is
$SL_{M+1}$-invariant, a map $\nu: \barr{P_{g,e}} \rarr \barr{U_g(e,1)}$ is
obtained.

It remains to prove $\nu$ is an isomorphism. Since $\barr{U_g(e,1)}$ is
normal by Proposition (\ref{popo}), it suffices to show $\nu$ is bijective.
Surjectivity is clear. Injectivity is established in section (\ref{ender}).

\subsection{The Construction of $\nu$}
Multidegree (semi)stability corresponds to
slope-(semi)stability in the following manner:
\begin{lm}
\label{ms}
Let $C$ be a Deligne-Mumford quasi-stable curve.
If $L$ is a very ample, degree $e$ line bundle on $C$ of
(semi)stable multidegree, then
$E=\psi_*(L)$ is a slope-(semi)stable, torsion free
sheaf of uniform rank $1$ and degree $e$ on $C_s$. Also, if
$L$ is of strictly semistable multidegree, then $E$ is strictly
slope-semistable.
\end{lm}
First we need a simple technical result. For each complete
subcurve $D$ of $C$,
define the sheaf $F_D$ on $C$  by the sequence:
$$ 0 \rarr F_D \rarr L \rarr L_{D^c} \rarr 0.$$
$F_D$ is the subsheaf of sections of $L$ with support on $D$.
In fact, $F_D$ is exactly the subsheaf of sections of $L_D$
vanishing on $D \cap D^c$. Therefore $degree(F_D)=e_D-k_D$.
Note by Riemann-Roch, $\chi(F_D)=degree(F_D)+1-g_D$ where
$g_D$ is the arithmetic genus of $D$. We obtain,
\begin{equation}
\label{dami}
e_D= \chi(F_D) + g_D -1 +k_D.
\end{equation}
\begin{lm}
\label{tec}
\sloppy
Let $C$ be a Deligne-Mumford quasi-stable curve.
Let $L$ be a very ample line bundle of semistable multidegree.
For every complete subcurve $D\subset C$, $R^1\psi_*(F_D)=0$.
\end{lm}
\begin{pf}
A fiber  of $\psi$ is either a point  or a destabilizing
$\proj^1$. By inequality (\ref{basinq}) and ampleness,
the restriction of $L$ to a destabilizing $\proj^1$ is
 $\oh_{\proj^{1}}(1)$. Let $P$ be a destabilizing
$\proj^1$ of $C$. There are five cases
\begin{enumerate}
\item [(i.)] $P\subset D^c$, $P\cap D =\emptyset$. Then $F_D|_P=0$.
\item [(ii.)] $P\subset D^c$, $P\cap D \neq \emptyset$.
 Then $F_D|_P$ is torsion.
\item [(iii.)] $P\subset D$, $P \cap D^c = \emptyset$. Then $F_D|_P=\oh_P(1)$.
\item [(iv.)] $P\subset D$, $\#|P \cap D^c|=1$. Then
$F_D|_P= \oh_P$.
\item [(v.)] $P\subset D$, $\# |P \cap D^c|=2$. Then $F_D|_P=\oh_P(-1)$.
\end{enumerate}
In each case, $h^1(P,F_D|_P)=0$.
The
vanishing of $h^1(P,F_D|_P)$
and the simple character of $\psi$ imply $R^1\psi_*(F_D)=0$.
\end{pf}

\begin{pf} [Of Lemma (\ref{ms})]
Let $C^1$ be the union of the destabilizing $\proj^1 \ 's$ of
$C$. By Lemma (\ref{tec}), $R^1\phi_*(F_{C^1})=0$.
Since each destabilizing
$P$ in $C^1$ is of type (v) in the proof of
Lemma (\ref{tec}),
$F_{C^1}|_P=\oh_{P}(-1)$. It follows that
$\psi_*(F_{C^1}) = 0.$
By the long exact sequence associated to
$$0 \rarr F_{C^1} \rarr L \rarr L_{C^{1c}} \rarr 0,$$
it follows
$E \cong \psi_*( L_{C^{1c}})$. Since $L_{C^{1c}}$ is torsion
free on $C^{1c}$ and the morphism $C^{1c} \rarr C_s$ is finite,
$E$ is torsion free. Certainly, E is of uniform rank $1$.

Let $0 \rarr G \rarr E$
be a proper subsheaf. Let $D_s\subset C_s$ be the support of
$G$. If $D_s=C_s$, the inequality of slope-stability follows
trivially. We can assume $D_s$ is a complete, proper subcurve. Let
$D= \psi^{-1}(D_s)$. $D\subset C$ is a complete, proper subcurve.
It is clear that
$G$ is a subsheaf of $\psi_*(F_D)$ with torsion quotient. Therefore, it
suffices to check slope-(semi)stability for $\psi_*(F_D)$. Certainly
$$h^0(C_s,\psi_*(F_D)) = h^0(C,F_D).$$
By Lemma (\ref{tec}), $R^1\psi_*(F_D)=0$. Thus by a
degenerate Leray spectral sequence ([H], Ex. 8.1, p.252),
$$h^1(C_s,\psi_*(F_D)) = h^1(C,F_D).$$
Hence $\chi(F_D)=\chi(\psi_*(F_D))$.
Similarly $\chi(L)=\chi(\psi_*(L))$.
Inequality (\ref{basinq}) and equation (\ref{dami}) yield:
$$(\chi(F_D)+g_D-1+k_D)   - (\chi(L)+g-1)
\cdot \paren{2g_D-2+k_D \over 2g-2}\ (\leq)
< \  k_D/2.$$
After some manipulation, we see
$${\chi(F_D) \over \om_{D,C}} \ \ (\leq) < \ \ {\chi(L)\over 2g-2}\ .$$
The above results yield
$${\chi(\psi_*(F_D)) \over \om_{D_s,C_s}}
 \ \ (\leq) < \ \ {\chi(\psi_*(L)) \over 2g-2}\ .$$
Hence, $\psi_*(L)$ is slope-(semi)stable. The final claim about
strict semistability also follows from the proof.
\end{pf}

Let $\psi: U_Z \rarr U_Z^s$, $L=\oh_{U_Z}(1)$, and $\cal{E}=\psi_*(L)$ be
as defined in section (\ref{fofo}). We now establish that $\cal{E}$ is flat
over
$Z_{g,e}$. The vanishing of $R^1\psi_*(L)$ is proved by
Zariski's theorem on formal functions in Lemma (\ref{zyz}).
The flatness criterion of Lemma (\ref{tpt}) is
then applied to obtain the required flatness.

First we need an auxiliary result.
Let $[C]\in Z_{e,g}$  and
let $P\subset U_Z$ be a destabilizing
$\proj^1$ of $C$. The
conormal bundle, $N_P^*$, of $P$ in $U_Z$  is locally free ($P$, $U_Z$ are
nonsingular).
Recall
a locally free sheaf $\bigoplus \oh_{\proj^1}(a_i)$ on
$\proj^1$ is said to be non-negative if each
$a_i\geq 0$.
\begin{lm}
\label{tyt}
$N_P^*$ is non-negative.
\end{lm}
\begin{pf}
Let $T_U$ and $T_Z$ denote the
tangent sheaves of $U_Z$ and $Z_{e,g}$. Let $\rho:U_Z\rarr Z_{e,g}$
be the natural morphism. There is  a differential map
$$d\rho: T_U \rarr \rho^*(T_Z).$$
Restriction to $P$ yields a (non-exact) sequence
$$0\rarr T_P \rarr T_U|_P \rarr \rho^*(T_Z)|_P .$$
Certainly, $\rho^*(T_Z)|_P \cong \bigoplus \oh_P$.
We obtain a map
$$\alpha: N_P\rarr \bigoplus \oh_P.$$
Let $\hat{P}\subset P$ be the locus of
nonsingular point of $C$.  Since the morphism $\rho$ is smooth on $\hat{P}$,
$\alpha |_{ \hat{P}}$ is an isomorphism of sheaves. Since
$N_P$ is a torsion free sheaf, $\alpha$ must be an
injection of sheaves. It follows easily $N_P$ is
non-positive. Hence $N_P^*$ is non-negative.
\end{pf}
In fact, an examination of the deformation theory yields
$N_P^* \cong \oh_P(1) \oplus \oh_P(1) \oplus I$
where $I$ is a trivial bundle. We will need only
the non-negativity result.

\begin{lm}
\label{zyz}
$R^1\psi_*(L)=0.$
\end{lm}
\begin{pf}
Let $\zeta \in U_{Z}^s$. It suffices to prove
$$R^1\psi_*(L)_\zeta=0$$
in case $\zeta$ is a node of stable curve destabilized
in $U_Z$. Let $m$ be the ideal of $\zeta$.
$\psi^{-1}(m)$ is the ideal of the nonsingular
destabilizing $P=\proj^1$.
Let $P_n$ denote the subscheme of $U_Z$  defined
by $\psi^{-1}(m^n) \cong \psi^{-1}(m)^n$. Let
$L_n$ be the restriction of $L$ to $P_n$.
By
Zariski's Theorem on formal functions:
$$R^1\psi_*(L)_\zeta\hat{\ } \ \ \cong \ \ \stackrel {lim}{\leftarrow}
H^1(P_n,L_n).$$
Since completion is faithfully flat for noetherian local
rings,
it suffices to show for each $n\geq 1$, $h^1(P_n,L_n)=0$. As above,
denote the conormal bundle of $P$ in $U_Z$ by
$N_P^*$. Since the varieties in question  are nonsingular, there
is an isomorphism on $P$:
$$m^{n-1}/m^n \cong Sym^{n-1} (N_P^*).$$

Since the pair $(C,L_C)$ is of semistable multidegree and
$P$ is a destabilizing $\proj^1$,
$L_1\cong \oh_P(1)$. Hence $h^1(P_1, L_1)=0.$
There is an exact sequence on $P_n$ for each $n\geq 2$:
$$0 \rarr m^{n-1}/m^n \otimes L_n \rarr L_n \rarr L_{n-1} \rarr 0.$$
There is a natural  identification
$$m^{n-1}/m^n\otimes L_n \cong Sym^{n-1}(N_P^*) \otimes_{\oh_P} \oh_P(1).$$
{}From the non-negativity of $N_P^*$ (Lemma (\ref{tyt})), we see
$$h^1(P_n,m^{n-1}/m^n \otimes L_n)=0.$$
 By the induction hypothesis $$h^1(P_n,L_{n-1})=h^1(P_{n-1},L_{n-1})=0.$$
Thus $h^1(P_n,L_n)=0$.
The proof is complete.
\end{pf}

\begin{lm}
\label{tpt}
Let $\phi: B_1 \rarr B_2$ be a projective morphism of schemes over $A$.
If $F$ is a sheaf on $B_1$ flat over $A$ and
$\forall \ i\geq 1$, $\ R^i\phi_*(F)=0$, then
$\phi_*(F)$ is flat over $A$.
\end{lm}
\begin{pf}
We can assume $A$ and $B_2$ are affine and $B_1 \cong \proj^k_{B_2}$.
Let $\cal{U}$ be the standard $k+1$ affine cover of $B_2$.
There is a Cech resolution computing the cohomology of
$F$ on $B_1$:
\begin{equation}
\label{cech}
0\rarr H^0(B_1,F) \rarr C^0(\cal{U},F) \rarr C^1(\cal{U},F) \rarr
\ldots \rarr C^k(\cal{U},F) \rarr 0.
\end{equation}
Since $R^i\phi_*(F)=0$ for $i\geq 1$, the resolution
(\ref{cech}) is exact. Since $F$ is $A$-flat, the
Cech modules $C^j(\cal{U},F)$ are all $A$-flat.  Exactness
of (\ref{cech}) implies $H^0(B_1,F) \cong \phi_*(F)$ is
$A$-flat.
\end{pf}

\begin{lm}
\label{flaty}
$\cal{E}$ is a flat family of slope-semistable, torsion free
sheaves of uniform rank $1$ over $Z_{e,g}$.
\end{lm}
\begin{pf}
By Lemmas (\ref{zyz}) and (\ref{tpt}), $\cal{E}$ is
flat over $Z_{e,g}$.
Let $[C]\in Z_{e,g}$. We have a diagram:
\begin{equation}
\begin{CD}
C @>{i_C}>> U_Z \\
@VV{\psi_C}V  @VV{\psi}V \\
C_s @>>{i_{C_s}}> U^s_Z
\end{CD}
\end{equation}
If $i_{C_s}^*(\cal{E}) \cong \psi_{C*}(i_C^*(L))$,
then the proof is complete  by Lemma (\ref{ms}).
There is a natural morphism of sheaves
$$\gamma_C : i_{C_s}^*(\cal{E}) \rarr
\psi_{C*}(i_C^*(L)). $$
We first show $\gamma_C$ is a surjection. Since, by Lemma (\ref{ms}),
$\psi_{C*}(i_C^*(L))$ is a slope-semistable torsion free
sheaf of degree $e$, $\psi_{C*}(i_C^*(L))$ is generated by
global sections. There is a natural identification
$$H^0(C_s,\psi_{C*}(i_C^*(L)))\cong H^0(C,\oh_C(1)).$$
By the nondegeneracy
of $C$ and the nonspeciality of $\oh_C(1)$,
$H^0(C,\oh_C(1))$ is canonically isomorphic
to $H^0(\proj^M,\oh_{\proj^M}(1))$.  These sections extend
over $U_Z$ and thus appear in $i_{C_s}^*(\cal{E})$.
Therefore, $\gamma_C$ is a surjection.

Since $\psi$ is an isomorphism except at
destabilizing $\proj^{1} \ 's$ of $U_Z$,
the kernel of $\gamma_C$ is a torsion sheaf on $C_s$.
Flatness of
$\cal{E}$ over $Z_{e,g}$ implies
$\chi(i_{C_s}^*(\cal{E}))$ is independent of $[C]\in Z_{e,g}$.
By Lemma (\ref{ms}), $\chi(\psi_{C*}(i_C^*(L)))$ is independent
of $[C]\in Z_{e,g}$.
Over the  open locus of nonsingular curves  $[C]\in Z_{e,g}$,
 $\psi$ is an isomorphism
thus:
\begin{equation}
\label{comb}
\chi(i_{C_s}^*(\cal{E}))  =
\chi(\psi_{C*}(i_C^*(L)) ).
\end{equation}
 By the above considerations,  (\ref{comb})
holds for every $[C] \in Z_{e,g}$. Hence, the torsion kernel
of $\gamma_C$ must be zero.
$\gamma_C$ is an isomorphism. The proof is complete.
\end{pf}
By combining Lemma (\ref{flaty}) with Theorem (\ref{thwee}),
there exists a natural  morphism $\nu_Z: Z_{e,g} \rarr
\barr{U_g(e,1)}$. Since $\nu_Z$ is $SL_{M+1}$-invariant, $\nu_Z$ descends to
$\nu: \barr{P_{g,e}} \rarr \barr{U_g(e,1)}$.
Certainly  $\nu$ is surjective.
Since $\barr{U_g(e,1)}$ is normal, $\nu$ is an isomorphism
if and only if $\nu$ is injective. The injectivity of $\nu$ will
be established in section (\ref{ender}).

\subsection{Injectivity of $\nu$}
\label{ender}
Let $C$ be a Deligne-Mumford quasi-stable genus $g$ curve.
Let $D\subset C$ be a complete subcurve.
A node $z\in D$ is an {\em external} node of $D$ if
$z \in D^c$. $P\subset D$ is a {\em destabilizing $\proj^1$ of
$D$} if $P$ is a destabilizing $\proj^1$ of $C$.
A  destabilizing
$\proj^1$ of $D$ is {\em external} if
$P \cap D^c\neq \emptyset$.
Let $\hat{D}$ denote   $D$ minus all the
external destabilizing $\proj^1 \ 's$ of $D$.

$(C,L)$ is a {\em semistable pair} if $C$ is Deligne-Mumford
quasi-stable and
$L$ is a very ample line bundle of semistable
multirank. The semistable  pairs
$(C,L)$ and $(C',L')$ are isomorphic if there exists an
isomorphism $\gamma:C\rarr C'$ such that $\gamma ^*(L')\cong L$.
A complete subcurve
 $D\subset C$
is an {\em extremal subcurve} of the semistable pair
$(C,L)$ if equality holds in (\ref{basinq}):
$$e_D   - e \cdot \paren{ \om_{C,D} \over  2g-2} \  = \  k_D/2.$$
  The semistable pair $(C,L)$ is said to be {\em maximal}
if the following condition is satisfied:
if $z \in C$ is an external node of an extremal subcurve,
$z$ is contained in a destabilizing $\proj^1$.
Let $\psi: C \rarr C_s$ be the stable contraction.
By Lemma (\ref{ms}), $\psi_*(L)$ is a slope-semistable,
torsion free sheaf of uniform rank $1$. Let
$J(\psi_*(L))$ be the associated set of slope-stable Jordan-Holder
factors.

$(C,J)$ is a {\em Jordan-Holder pair} if $C$ is a Deligne-Mumford
stable curve and $J$ is a set of slope-stable, torsion free
sheaves. As before, the Jordan-Holder pairs $(C,J)$ and $(C',J')$
are isomorphic if there exists  an isomorphism $\gamma:C \rarr C'$ such that
$\gamma ^*(J') \cong J$.

\begin{lm}
\label{det}
Let $(C,L)$, $(C',L')$  be maximal semistable pairs. If
$(C_s, J(\psi_*(L)))$ and $(C'_s, J(\psi'_*(L')))$ are
isomorphic Jordan-Holder pairs,  then $(C,L)$ and $(C',L')$
are isomorphic semistable pairs.
\end{lm}
\begin{pf}
Consider a Jordan-Holder filtration of $E=\psi_*(L)$ on $C_s$:
$$0 = E_0 \subset E_1 \subset E_2  \subset \ldots \subset E_n=E.$$
Let $A_i=Supp(E_i)$. For each $1\leq i \leq n$,
$${\chi(E_i) \over w_{A_i,C_s}} = {\chi(E) \over 2g-2}.$$
 By the proof of Lemma (\ref{ms}), we
see the $B_i=\psi^{-1}(A_i)$ are extremal subcurves of $(C,L)$
for $1\leq i\leq n$.
As before,
let $F_{B_i}$ be the subsheaf of sections of $L$ with support
on $B_i$. From the proof of Lemma (\ref{ms}),
it follows $E_i \cong \psi_*(F_{B_i})$. For $1\leq i \leq n$,
let $X_i= \barr{A_i \setminus A_{i-1}}$ and
 $Y_i=\barr{B_i \setminus B_{i-1}}$. We see
$Supp(E_i/E_{i-1})=X_i$. If $z\in X_i$ is an internal node of $X_i$,
$z$ is destabilized by $\psi$ if and only if $(E_i/E_{i-1})$ is
locally isomorphic to $m_z$ at $z$. If $z\in X_i$ is an
external node, then there are two cases.
If $z\in A_{i-1}$, then $\psi^{-1}(z) \cong \bold{P}^1 \not\subset Y_i.$
If $z\in A_i^c$, then $\psi^{-1}(z) \cong  \bold{P}^1 \subset Y_i.$
These conclusions follow from the maximality of $(C,L)$.
It is now clear  $C_s$ and the Jordan-Holder factor $E_i/E_{i-1}$ determine
$\hat{Y_i}$ completely. Also, the $\hat{Y_i}$ are connected by
destabilizing $\proj^1 \ 's$. We have shown
$(C_s, J)$ determines $C$ up to isomorphism.
We will show
below in Lemmas (\ref{mjk}-\ref{pushiso}) that $L_{\hat{Y_i}}$
is determined up to isomorphism by $E_i/E_{i-1}$.
Since the $\hat{Y}_i$ are connected by destabilizing ${\proj^{1}}\ 's$, the
isomorphism class of the pair $(C,L)$ is determined by the
line bundles
$L_{\hat{Y_i}}$.
This completes the proof of the Lemma.
\end{pf}

\begin{lm}
\label{mjk}
There is an isomorphism
$\psi_*(L_{\hat{Y_i}}) \cong E_i/E_{i-1}.$
\end{lm}
\begin{pf}
We keep the notation of the previous Lemma.
Certainly $Supp(F_{B_i}/F_{B_{i-1}})=Y_i$.
Let $p_{i}$ be the divisor $B_{i-1}\cap Y_i \subset Y_i$.
Let $q_i$ be the divisor $ B_i^c\cap Y_i \subset Y_i$.
Since the points of $p_{i}$ lie on destabilizing $\proj^1$'s joining
$Y_i$ to $Y_{i-1}$ and the points of $q_{i}$ lie on
destabilizing $\proj^1$'s joining $Y_i$ to $Y_{i+1}$, we
note $p_i \cap q_i = \emptyset$.
There is a isomorphism $L_{Y_i} \cong F_{Y_i}\otimes \oh_{Y_i} (p_i +q_i)$.
Since there is  an exact sequence:
$$0 \rarr F_{Y_i} \rarr (F_{B_i}/F_{B_{i-1}}) \rarr p_i \rarr 0,$$
we see $F_{Y_i}= (F_{B_i}/F_{B_{i-1}})\otimes \oh_{Y_i}(-p_i)$.
Thus $$L_{Y_i} \cong (F_{B_i}/F_{B_{i-1}})\otimes \oh_{Y_i} (q_i).$$
Since $B_i$ is an extremal subcurve of $C$ and the pair $(C,L)$ is
maximal, we see $q_i$ lies on external destabilizing $\proj^1 \ 's$
of $Y_{i}$. Hence  $L_{\hat{Y_i}}$ is
isomorphic to
 $(F_{B_i}/F_{B_{i-1}})_{\hat{Y_i}}$.
We have an exact sequence on $Y_i$:
$$0 \rarr I_{\hat{Y_i}} \otimes (F_{B_i}/F_{B_{i-1}}) \rarr
(F_{B_i}/F_{B_{i-1}}) \rarr (F_{B_i}/F_{B_{i-1}})_{\hat{Y_i}} \rarr 0.$$
Let $P$ be an external destabilizing $\proj^1$ of $Y_i$.
If $P$ meets $B_{i-1}$ then $P\subset B_{i-1}$. Thus
each such $P$ meets $B_i^c$. It is now not hard to see
$ I_{\hat{Y_i}} \otimes (F_{B_i}/F_{B_{i-1}})$ restricts
to $\oh_P(-1)$ on each such $P$. Therefore, by familiar arguments,
$$\psi_*( (F_{B_i}/F_{B_{i-1}})_{\hat{Y_i}}) \cong
\psi_*( F_{B_i}/F_{B_{i-1}} ).$$
By Lemma (\ref{tec}),
$R^1\psi_*(F_{B_{i-1}})=0$. Hence
$$\psi_*(F_{B_i}/F_{B_{i-1}})
 \cong (\psi_*(F_{B_i})/\psi_*(F_{B_{i-1}})) \cong
E_i/E_{i-1}.$$
Following all the isomorphisms yields the Lemma.
\end{pf}

\begin{lm}
\label{pushiso}
Let $(C,L)$ be a semistable pair.
Let $D\subset C$ satisfying $\hat{D}=D$.
Then $L_D$ is determined up to isomorphism by $\psi_*(L_D)$.
\end{lm}
\begin{pf}
Let $D^1$ be  the union of the
destabilizing $\proj^1\ 's$ of $D$.
Note all these $\proj^1 \ 's$ are internal. Let
 $D'= \barr{D \setminus D^1}$ and
denote the restriction of $\psi$ to $D'$ by $\psi'$.
Consider the sequence on $D$:
$$0 \rarr I_{D'}\otimes L_D \rarr L_D \rarr L_{D'} \rarr 0.$$
Since $I_{D'}\otimes L_D$ restricts to $\oh_{\proj^1}(-1)$ on
each
destabilizing $\proj^1$ of $D$, we see
$$ \psi_* (I_{D'}\otimes L_D)= R^1 \psi_* (I_{D'} \otimes L_D)=0.$$
Thus $\psi'_*(L_{D'})\cong
 \psi_*(L_{D'}) \cong \psi_*(L_D)$. Since
$\psi'$ is a finite affine morphism,
$$\beta: \psi'\ ^*(\psi'_*(L_{D'})) \rarr L_{D'} \rarr 0.$$
Let $\tau $ be the torsion subsheaf of $\psi'\ ^*(\psi'_*(L_{D'}))$.
Since $\beta$ is generically an isomorphism and $L_{D'}$ is
torsion free on $D'$, we see
$$L_{D'} \cong  (\psi'\ ^*(\psi'_*(L_{D'}))/ \tau).$$
We have shown that $L_{D'}$ is determined up to isomorphism by
$\psi_*(L_D)$. It is clear that $L_{D'}$
determines $L_D$ up to isomorphism.
\end{pf}

Let $\rho: Z_{e,g} \rarr \barr {P_{g,e}}$ be the
quotient map.
Let $\zeta \in \barr{P_{g,e}}$. It follows from
the results of [Ca] (Lemma 6.1, p.640)
that there exists a $[C]\in Z_{g,e}$ satisfying:
\begin{enumerate}
\item[(i.)]$\rho([C]) = \zeta$.
\item[(ii.)] $(C, L=\oh_C(1))$ is a maximal semistable pair.
\end{enumerate}
Let $\psi_C:C\rarr C_s$ be the stable contraction.
Let $E=\psi_*(L)$.
Let $J$ be the Jordan-Holder factors of $E$ on $C_s$. From the
definition of $\nu$,
$\nu(\zeta)$ is the element of $\barr{U_g(e,1)}$ corresponding to
the isomorphism class of the data $(C_s,J)$.
 By Lemmas (\ref{mjk}-\ref{pushiso}), the isomorphism
class of $(C,L)$ is determined by the isomorphism class of
$(C_s, J)$. Therefore
$\nu$ is injective. By the previous discussion,
$\nu$ is an isomorphism.
\begin{tm} There is a natural isomorphism
$\nu: \barr{P_{g,e}} \rarr \barr{U_g(e,1)}$.
\end{tm}

\hbadness=9999

\noindent
\address{Department of Mathematics \\ University of Chicago \\
5734 S. University Ave. \\ Chicago, IL 60637 \\
rahul@@math.uchicago.edu}

\end{document}